\documentclass[11pt]{article}

\usepackage[a4paper,margin=1in]{geometry}

\usepackage{amsmath, amssymb, amsfonts, amsthm, mathtools}
\usepackage{graphicx}
\usepackage{algorithm}
\usepackage{algpseudocode}
\usepackage{subcaption}
\usepackage{mathrsfs}
\usepackage{bbm}
\usepackage{bm}
\usepackage{xcolor}
\usepackage{booktabs}
\usepackage{url}

\newtheorem{theorem}{Theorem}
\newtheorem{corollary}{Corollary}
\newtheorem{proposition}{Proposition}
\newtheorem{assumption}{Assumption}
\newtheorem{lemma}{Lemma}
\newtheorem{definition}{Definition}

\newcommand{\keywords}[1]{%
  \vspace{0.5em}
  \noindent\textbf{Keywords:} #1
}

\title{Law-Dependent Poisson-Refresh Randomised Feedback Control for Stochastic Systems}

\author{
Randa Herzallah\\
Warwick Mathematics Institute, University of Warwick, \\ Coventry CV4 7AL, United Kingdom
}

\date{}

\begin{document}

\maketitle

\begin{abstract}
This paper studies a continuous-time probabilistic control problem for
stochastic systems in which the feedback law assigns a probability distribution
over control inputs. The control is modelled as a law-dependent Poisson-refresh
process: it is held constant between refresh times and resampled at Poisson
times from a feedback rule depending on the current state and state law. The
closed-loop system is therefore formulated as a jump--diffusion for the joint
state--control process.
The objective is to design the sampling rule by comparing the realised
closed-loop system with a prescribed reference description. The cost separates
two local terms: a relative-entropy penalty for the instantaneous refresh kernel
and a state-side drift-mismatch cost. The state-side term is evaluated at the
realised control value and is motivated by a frozen-control Girsanov identity.
Under global Lipschitz, moment, and uniform ellipticity assumptions, we prove
well-posedness of the closed loop and a dynamic programming principle on the
joint state--control law space.
For linear data and admissible mean-dependent affine-Gaussian refresh rules,
the full joint law need not remain Gaussian, but the first and second moments
of the joint state--control process close exactly. This yields a deterministic
reduced control problem on moment space, with an explicit running cost, a
reduced Hamiltonian, and a verification theorem.
\end{abstract}


\section{Introduction}

A large body of work has studied feedback control for stochastic systems in
both continuous and discrete time. In the classical formulation, the objective
is usually to choose a control policy that minimises an expected cost functional
\cite{Charalambous95}. For linear systems with Gaussian noise and quadratic
costs, this leads to linear quadratic Gaussian theory
\cite{HUANG201468,ZHANG2023105551,SUN2023100899,SUN2016,Wonham68}. This gives
a highly tractable benchmark, but outside the linear--Gaussian setting explicit
solutions are rarely available, and the associated dynamic programming
equations are typically nonlinear and difficult to solve.

A related but different line of work formulates stochastic control in
information-theoretic terms. Path-integral control, Kullback--Leibler control,
entropy-regularised control, and Schr\"odinger bridge formulations all use
relative entropy or exponential reweighting ideas to compare controlled and
reference probabilistic evolutions
\cite{YANG20146994,KAZIM2024100931,Theodorou12,Theodorou2013FromIT,Ito22,
ChenGeorgiouPavon2021,JiaZhou2022TD,JiaZhou2023QLearning}. These approaches
show that, under suitable structure, stochastic control problems can be related
to changes of measure, free-energy identities, or controlled transport between
probability laws. However, many such formulations either compare full path laws,
work with transition-kernel modifications, or introduce entropy as a regulariser
on top of an otherwise standard control objective. They do not directly model a
continuous-time controller as a sampled control process that is physically held
between random refresh times.

Another important development is the study of law-dependent and
McKean--Vlasov control. In this setting, the state law itself becomes part of
the continuation variable, and dynamic programming is naturally formulated on
spaces of probability measures
\cite{Annunziato18,Annunziato2010,Cardaliaguet2019,Carmona2019Modeling,
Bahlali2017RelaxedMF,DjetePossamaiTan2022MKVDPP,PhamWei2017}. This provides a
powerful framework for mean-field interactions and distribution-dependent
dynamics. More recently, random-measure approaches to continuous-time
reinforcement learning have given a rigorous interpretation of sampled relaxed
controls through grid sampling and limiting random-measure SDEs
\cite{BenderThuan2026RMARL}. These approaches are primarily concerned
with the analysis, learning, or execution of randomised or relaxed controls,
rather than with modelling the realised sampled controller itself as part of the
closed-loop dynamics.

The present paper addresses this gap by studying a continuous-time stochastic
control problem in which the randomised feedback controller is realised as an
actual sampled process. The control \(U\) is held constant between refresh times
and is resampled at Poisson times from a law-dependent feedback rule. More
precisely, at a refresh time, the new control value is sampled from a kernel
depending on the current state and the current state law. The closed-loop system
is therefore formulated as a jump--diffusion for the joint process \((X,U)\):
the state evolves according to the realised value of the control, while the
control component evolves through Poisson refreshes.

This differs from standard entropy-regularised and path-space formulations in
two main ways. First, the primitive control object is not only a probability
kernel, but a sampling rule that generates a piecewise-constant control process.
Second, the natural continuation variable is the joint law of state and realised
control, \(\mathcal L(X_t,U_t)\) rather than the state law alone. This is important because the state-side cost depends on the realised control value, not merely on the post-refresh kernel.

To optimise the sampling rule, we assign a cost to how the realised closed loop
deviates from a prescribed reference description. The objective functional is
built from two local terms. The first is a control-side relative-entropy cost
for the instantaneous refresh kernel relative to a prescribed reference refresh
rule. This term measures the deviation of the post-refresh control law from the
reference law and is independent of the particular parametrisation used to
sample from that law. The second is a state-side drift-mismatch cost induced by
a prescribed reference drift. This term is evaluated at the realised control
value and is motivated by a frozen-control Girsanov identity. Together, these
two terms give an additive finite-horizon cost for the realised closed-loop
system.

The first contribution of the paper is to place this Poisson-refresh sampled
formulation on a precise mathematical foundation. Under global Lipschitz,
moment, and uniform ellipticity assumptions, we introduce an admissible class of
sampling rules and prove closed-loop well-posedness for the joint
state--control process. We also prove stability with respect to the initial
joint law, continuity of the induced state-law flow, c\`adl\`ag regularity of
the joint law flow, and stability of admissibility under restriction and
concatenation.

The second contribution is the formulation of the corresponding continuation
problem on
\(
  [0,T]\times\mathcal P_2(\mathbb R^n\times\mathbb R^q).
\)
Using the restriction and concatenation properties of the admissible class, we
prove a dynamic programming principle for the value function with the joint law
of \((X_t,U_t)\) as continuation variable.

The third contribution is the identification of the two local information
costs. We define the control-side entropy at the level of instantaneous refresh
kernels and record an auxiliary marked-refresh representation of this entropy.
We also prove a frozen-control Girsanov identity that motivates the state-side
drift-mismatch term. This separates the cost of changing the refresh mechanism
from the cost of changing the induced state dynamics.

Finally, we study a linear-data example with mean-dependent affine-Gaussian
refresh rules. In the exact Poisson-refresh setting, the full joint law of
\((X_t,U_t)\) is not generally Gaussian, because the control is reset at random
times and then held between refreshes. Nevertheless, within the admissible
mean-dependent affine-Gaussian class, the first and second moments of
\((X_t,U_t)\) close exactly. This yields a deterministic reduced control problem
on moment space, with an explicit running cost, a reduced dynamic programming
principle, a reduced Hamiltonian, and a verification theorem.

The main contributions of the paper are summarised as follows:
\begin{enumerate}
\item
A Poisson-refresh formulation of law-dependent randomised feedback in which the
control is an actual sampled process and the continuation variable is the joint
state--control law.

\item
A well-posedness and dynamic programming theory for the resulting closed-loop
problem under a concrete admissible class of sampling rules.

\item 
A finite-horizon information-cost objective combining refresh-kernel entropy
with a state-side drift-mismatch cost motivated by a frozen-control Girsanov
identity.

\item
An exact moment reduction for linear data under admissible mean-dependent
affine-Gaussian refresh rules, leading to a deterministic reduced control
problem on moment space.
\end{enumerate}

The paper is organised as follows.
Section~\ref{sec:poisson-refresh-formulation} introduces the Poisson-refresh
randomised feedback formulation and the closed-loop jump--diffusion model.
Section~\ref{sec:control-entropy} defines the local control-side entropy cost.
Section~\ref{sec:state-cost} introduces the reference drift, the state-side
cost, and the finite-horizon objective. Section~\ref{sec:dpp-poisson}
formulates the continuation problem and proves the dynamic programming
principle. Section~\ref{sec:LG-moment} studies the linear-data moment reduction
under mean-dependent affine-Gaussian refresh rules. Concluding remarks are given
in Section~\ref{sec:conc_Rem}.

\section{Poisson-refresh randomized feedback controls and closed-loop jump--diffusion formulation}
\label{sec:poisson-refresh-formulation}

We consider a finite-horizon control problem on the interval \([0,T]\). The
state process \(X\) takes values in \(\mathbb R^n\), the control process \(U\)
takes values in \(\mathbb R^q\), and the state dynamics are governed by the stochastic
differential equation (SDE)
\begin{equation}
  \mathrm dX_t
  =
  f(X_t,U_t)\,\mathrm dt
  +
  \sigma(X_t,U_t)\,\mathrm dW_t,
  \label{eq:primitive-controlled-SDE-poisson}
\end{equation}
where \(W\) is an \(n\)-dimensional standard Wiener process,
\(f:\mathbb R^n\times\mathbb R^q\to\mathbb R^n\), and
\(\sigma:\mathbb R^n\times\mathbb R^q\to\mathbb R^{n\times n}\).

The control variable in this paper is modelled as a sampled process
\(U\), piecewise constant between refresh times and refreshed at Poisson times
from a law-dependent randomized feedback rule. 

For \(d\in\mathbb N\), let \(\mathcal P(\mathbb R^d)\) denote the set of Borel
probability measures on \(\mathbb R^d\), and let
\[
  m_2(\nu):=\int_{\mathbb R^d}\|y\|^2\,\nu(\mathrm dy),
  \qquad
  \mathcal P_2(\mathbb R^d)
  :=
  \bigl\{
    \nu\in\mathcal P(\mathbb R^d): m_2(\nu)<\infty
  \bigr\}.
\]
We write \(W_2\) for the \(2\)-Wasserstein distance on
\(\mathcal P_2(\mathbb R^d)\).

For \(0\le s_0\le s_1\le T\), let
\(\mathcal X_{s_0,s_1}:=C([s_0,s_1];\mathbb R^n)\) and
\(\mathcal U_{s_0,s_1}:=D([s_0,s_1];\mathbb R^q)\),
and equip both spaces with their standard Borel \(\sigma\)-fields. Set
\(\mathcal Y_{s_0,s_1}:=\mathcal X_{s_0,s_1}\times\mathcal U_{s_0,s_1}\).
In particular, we write
\(\mathcal X:=\mathcal X_{0,T}\),
\(\mathcal U:=\mathcal U_{0,T}\), and
\(\mathcal Y:=\mathcal Y_{0,T}\).
Whenever \(s_0\le s_1\le s_2\), we denote by
\(\Pi_{s_0,s_1}:\mathcal Y_{s_0,s_2}\to\mathcal Y_{s_0,s_1}\)
the canonical restriction map
\(\Pi_{s_0,s_1}(x,u):=(x|_{[s_0,s_1]},u|_{[s_0,s_1]})\).
We also set
\(\Sigma(x,u):=\sigma(x,u)\sigma(x,u)^\top\) for
\((x,u)\in\mathbb R^n\times\mathbb R^q\).

Throughout the paper we assume the following. 
\begin{assumption}\label{assump:primitive-coeff}
The following hold.
\begin{enumerate}
\item
The functions \(f\) and \(\sigma\) are globally Lipschitz in \((x,u)\).

\item
There exist constants \(0<\underline{c}_\Sigma\le \overline{c}_\Sigma<\infty\) such
that, for every \((x,u)\in\mathbb R^n\times\mathbb R^q\),
\begin{equation}
  \underline{c}_\Sigma I_n
  \le
  \Sigma(x,u)
  \le
  \overline{c}_\Sigma I_n
  \label{eq:uniform-ellipticity-poisson}
\end{equation}
in the sense of symmetric matrices.

\item
A refresh rate \(\lambda>0\) is fixed once and for all.
\end{enumerate}
\end{assumption}

In the present Poisson-refresh formulation, the admissible control object is a
sampling rule, denoted by \(\alpha\). This rule prescribes how the
control process is refreshed at Poisson times. More precisely, at a refresh
time \(s\), the new control value is generated as
\(
  \alpha_s(x,\mu,z),
\)
where \(x\in\mathbb R^n\) is the current state,
\(\mu\in\mathcal P_2(\mathbb R^n)\) is the current state law, and
\(z\in[0,1]\) is a fresh auxiliary randomization variable. Between refresh
times, the control process \(U\) remains constant. Thus \(\alpha\) is the
control variable of the optimization problem, while \(U\) is the actual sampled
control process generated from \(\alpha\).

Fix constants \(L_\alpha,M_\alpha>0\). For \(0\le t\le T\), let \(\mathfrak A[t,T]\) denote the class of maps
\(\alpha=(\alpha_s)_{s\in[t,T]}\), with
\(
  \alpha_s:\mathbb R^n\times\mathcal P_2(\mathbb R^n)\times[0,1]
  \to\mathbb R^q,
\)
such that the map \((s,x,\mu,z)\mapsto\alpha_s(x,\mu,z)\) is Borel measurable
on \(
  [t,T]\times \mathbb R^n\times \mathcal P_2(\mathbb R^n)\times [0,1]
\), and the following estimates hold for all \(s\in[t,T]\), all
\(x,x'\in\mathbb R^n\), and all \(\mu,\mu'\in\mathcal P_2(\mathbb R^n)\):
\begin{align}
  \int_0^1
    \|\alpha_s(x,\mu,z)-\alpha_s(x',\mu',z)\|^2
  \,\mathrm dz
  &\le
  L_\alpha^2
  \Bigl(
    \|x-x'\|+W_2(\mu,\mu')
  \Bigr)^2,
  \label{eq:sampling-rule-lipschitz}
  \\
  \int_0^1
    \|\alpha_s(x,\mu,z)\|^2
  \,\mathrm dz
  &\le
  M_\alpha
  \Bigl(
    1+\|x\|^2+m_2(\mu)
  \Bigr).
  \label{eq:sampling-rule-second-moment}
\end{align}
We also write \(\mathfrak A:=\mathfrak A[0,T]\).

Every admissible sampling rule \(\alpha\in\mathfrak A[t,T]\) induces an
instantaneous feedback kernel by pushing forward Lebesgue measure on \([0,1]\)
through the map \(z\mapsto \alpha_s(x,\mu,z)\). Equivalently, for each
\((s,x,\mu)\), it specifies the conditional law of the post-refresh control
value at time \(s\). Thus the law of the refreshed control is the probability
measure on \(\mathbb R^q\) generated by the sampling rule at the point
\((s,x,\mu)\). 

The condition in \eqref{eq:sampling-rule-lipschitz} is an
\(L^2([0,1])\)-Lipschitz condition in the auxiliary sampling variable. This is
the form needed later when estimating the Poisson jump coefficient
\(\mathfrak j_s(x,u,\mu,z):=\alpha_s(x,\mu,z)-u\). The same \(L^2\)-structure also gives useful estimates for the probability
kernel induced by the sampling rule, as recorded in the following proposition.

\begin{proposition}\label{prop:kernel-induced-by-sampling-rule}
For \(\alpha\in\mathfrak A[t,T]\), define
\begin{equation}
  c^\alpha_s(B\mid x,\mu)
  :=
  \int_0^1
    \mathbf 1_B\!\bigl(\alpha_s(x,\mu,z)\bigr)\,\mathrm dz,
  \qquad
  B\in\mathcal B(\mathbb R^q).
  \label{eq:kernel-induced-by-sampling-rule}
\end{equation}
Then \(c^\alpha=(c^\alpha_s)_{s\in[t,T]}\) is a Borel probability-kernel family
on \(\mathbb R^q\), and for all \(s\in[t,T]\), \(x,x'\in\mathbb R^n\), and
\(\mu,\mu'\in\mathcal P_2(\mathbb R^n)\),
\begin{align}
  W_2\!\bigl(
    c^\alpha_s(\cdot\mid x,\mu),
    c^\alpha_s(\cdot\mid x',\mu')
  \bigr)
  &\le
  L_\alpha
  \Bigl(
    \|x-x'\|+W_2(\mu,\mu')
  \Bigr),
  \label{eq:kernel-from-alpha-lipschitz}
  \\
  \int_{\mathbb R^q}\|u\|^2\,c^\alpha_s(\mathrm du\mid x,\mu)
  &\le
  M_\alpha
  \Bigl(
    1+\|x\|^2+m_2(\mu)
  \Bigr).
  \label{eq:kernel-from-alpha-moment}
\end{align}
\end{proposition}

\begin{proof}
For each Borel set \(B\subseteq\mathbb R^q\), the map
\(
  (s,x,\mu)\mapsto c_s^\alpha(B\mid x,\mu)
  =
  \int_0^1 \mathbf 1_B\!\bigl(\alpha_s(x,\mu,z)\bigr)\,\mathrm dz
\)
is Borel measurable by measurability of \(\alpha\) and Fubini's theorem \cite{Kallenberg2021FMP,BertsekasShreve1978SOC}. Hence
\(c^\alpha\) is a Borel probability-kernel family.

For \eqref{eq:kernel-from-alpha-lipschitz}, use the common parameter
\(z\in[0,1]\) to couple
\(c^\alpha_s(\cdot\mid x,\mu)\) and
\(c^\alpha_s(\cdot\mid x',\mu')\). Then
\[
  W_2^2\!\bigl(
    c^\alpha_s(\cdot\mid x,\mu),
    c^\alpha_s(\cdot\mid x',\mu')
  \bigr)
  \le
  \int_0^1
    \|\alpha_s(x,\mu,z)-\alpha_s(x',\mu',z)\|^2\,\mathrm dz,
\]
and \eqref{eq:kernel-from-alpha-lipschitz} follows from
\eqref{eq:sampling-rule-lipschitz}. Likewise,
\eqref{eq:kernel-from-alpha-moment} follows directly from
\eqref{eq:sampling-rule-second-moment}.
\end{proof}

We now specify the Poisson-refresh mechanism in closed-loop form. The model is
formulated directly in terms of the sampling rule \(\alpha\), since \(\alpha\)
generates the actual sampled control process. The associated kernel family
\(c^\alpha\) will be used later only when defining local information costs.

To construct the refresh mechanism, we use a Poisson random measure on
\([t,T]\times[0,1]\): its time component gives the refresh times, and its mark
\(z\in[0,1]\) selects the post-refresh control value through \(\alpha\). This leads to the following notion of weak closed-loop solution.

\begin{definition}\label{def:weak-closed-loop-solution-poisson}
Fix \(0\le t\le T\), \(\zeta\in\mathcal P_2(\mathbb R^n\times\mathbb R^q)\),
and \(\alpha\in\mathfrak A[t,T]\). A weak closed-loop solution for
\((t,\zeta,\alpha)\) consists of:
\begin{enumerate}
\item
a filtered probability space
\(
  (\Omega,\mathcal F,(\mathcal F_s)_{s\in[t,T]},\mathbb Q)
\)
satisfying the usual conditions;

\item
an \(n\)-dimensional standard Wiener process
\(W=(W_s)_{s\in[t,T]}\), Brownian with respect to
\((\mathcal F_s)_{s\in[t,T]}\);

\item
a Poisson random measure \(N(\mathrm ds,\mathrm dz)\) on
\([t,T]\times[0,1]\), adapted to \((\mathcal F_s)_{s\in[t,T]}\), with
intensity \(\lambda\,\mathrm ds\,\mathrm dz\), independent of \(W\), and with
independent increments relative to the filtration \((\mathcal F_s)_{s\in[t,T]}\);

\item
an \(\mathbb R^n\times\mathbb R^q\)-valued c\`{a}dl\`{a}g adapted process
\(
  (X,U)=((X_s,U_s))_{s\in[t,T]}
\)
such that \(X\) is continuous, \((X_t,U_t)\sim\zeta\),
\((X_t,U_t)\) is independent of
\(
  (W_s-W_t,\;N((t,s]\times B))_{s\in[t,T],\,B\in\mathcal B([0,1])}
\), and which also satisfies
\[
  \mathbb E^{\mathbb Q}\!\left[
    \sup_{s\in[t,T]}
    \bigl(
      \|X_s\|^2+\|U_s\|^2
    \bigr)
  \right]
  <\infty.
\]
Writing \(\mu_s:=\mathcal L_{\mathbb Q}(X_s)\) for \(s\in[t,T]\), one has
\begin{align}
  \mathrm dX_s
  &=
  f(X_s,U_s)\,\mathrm ds
  +
  \sigma(X_s,U_s)\,\mathrm dW_s,
  \label{eq:state-equation-poisson-refresh}
  \\
  U_s
  &=
  U_t
  +
  \int_{(t,s]\times[0,1]}
    \bigl(
      \alpha_r(X_{r-},\mu_r,z)-U_{r-}
    \bigr)
    N(\mathrm dr,\mathrm dz),
  \qquad
  s\in[t,T],
  \label{eq:control-refresh-equation}
\end{align}
\(\mathbb Q\)-almost surely.
\end{enumerate}

For any weak closed-loop solution
\(
  (\Omega,\mathcal F,(\mathcal F_s),\mathbb Q,W,N,X,U)
\)
for \((t,\zeta,\alpha)\), its path law on \(\mathcal Y_{t,T}\) is
\(
  \mathbb P:=\mathcal L_{\mathbb Q}(X,U)\in\mathcal P(\mathcal Y_{t,T}),
\)
and its associated joint law flow is
\(
  \eta_s:=\mathcal L_{\mathbb Q}(X_s,U_s),\)
 where \( s\in[t,T].
\)
The associated state-law flow is the family
\(
  (\mu_s)_{s\in[t,T]}
\)
introduced above.
\end{definition}

Equation \eqref{eq:control-refresh-equation} means that \(U\) is piecewise
constant between refresh times and that, at each jump time \(r\) of the Poisson
clock, the post-jump value \(U_r\) is reset to
\(\alpha_r(X_{r-},\mu_r,Z_r)\), where \(Z_r\in[0,1]\) is the corresponding
Poisson mark. The next theorem gives the well-posedness and basic stability estimates for the closed-loop Poisson-refresh dynamics.

\begin{theorem}\label{thm:wellposedness-poisson-refresh}
Under Assumption~\ref{assump:primitive-coeff}, for every \(0\le t\le T\), every
\(\zeta\in\mathcal P_2(\mathbb R^n\times\mathbb R^q)\), and every
\(\alpha\in\mathfrak A[t,T]\), there exists a weak closed-loop solution in the
sense of Definition~\ref{def:weak-closed-loop-solution-poisson} for
\((t,\zeta,\alpha)\). Moreover, the law of the pair \((X,U)\) on
\(\mathcal Y_{t,T}\) is unique. We denote this unique path law by
\(
  \mathbb P^{t,\zeta,\alpha}\in\mathcal P(\mathcal Y_{t,T}).
\)
For \(s\in[t,T]\), define
\(
  \eta_s^{t,\zeta,\alpha}
  :=
  \mathcal L_{\mathbb P^{t,\zeta,\alpha}}(X_s,U_s),\) and \(
  \mu_s^{t,\zeta,\alpha}
  :=
  \mathcal L_{\mathbb P^{t,\zeta,\alpha}}(X_s).
\)
Additionally, there exists a constant \(C>0\), depending only on the constants in
Assumption~\ref{assump:primitive-coeff}, on \(L_\alpha\), \(M_\alpha\), and on
\(T\), such that for every \(0\le t\le T\), every
\(\alpha\in\mathfrak A[t,T]\), every
\(\zeta,\zeta'\in\mathcal P_2(\mathbb R^n\times\mathbb R^q)\), and every
\(s,r\in[t,T]\),
\begin{align}
  \mathbb E^{t,\zeta,\alpha}\!\left[
    \sup_{\ell\in[t,T]}
    \bigl(
      \|X_\ell\|^2+\|U_\ell\|^2
    \bigr)
  \right]
  &\le
  C\bigl(1+m_2(\zeta)\bigr),
  \label{eq:moment-bound-joint-poisson}
  \\
  \sup_{\ell\in[t,T]}
  W_2\!\bigl(
    \eta_\ell^{t,\zeta,\alpha},
    \eta_\ell^{t,\zeta',\alpha}
  \bigr)
  &\le
  C\,W_2(\zeta,\zeta'),
  \label{eq:stability-in-initial-joint-law}
  \\
  W_2\!\bigl(
    \mu_s^{t,\zeta,\alpha},
    \mu_r^{t,\zeta,\alpha}
  \bigr)
  &\le
  C\bigl(1+m_2(\zeta)^{1/2}\bigr)|s-r|^{1/2}.
  \label{eq:holder-state-law-flow-poisson}
\end{align}
In particular,
\[
  s\longmapsto \mu_s^{t,\zeta,\alpha}
  \in
  C\bigl([t,T];\mathcal P_2(\mathbb R^n)\bigr),
\]
and
\[
  s\longmapsto \eta_s^{t,\zeta,\alpha}
  \in
  C\bigl([t,T];\mathcal P_2(\mathbb R^n\times\mathbb R^q)\bigr).
\]
\end{theorem}

\begin{proof}
Fix \(t\in[0,T]\), \(\zeta\in\mathcal P_2(\mathbb R^n\times\mathbb R^q)\), and
\(\alpha\in\mathfrak A[t,T]\). All stochastic constructions below are made on
filtered probability spaces satisfying the usual conditions. The Brownian
motion \(W\) and the Poisson random measure \(N\) have independent increments
relative to the filtration, and the initial random vector is independent of the
future increments of \((W,N)\).

Let
\[
  \mathcal M_t
  :=
  C\bigl([t,T];\mathcal P_2(\mathbb R^n)\bigr)
\]
equipped with
\[
  d_t(\varrho,\varrho')
  :=
  \sup_{s\in[t,T]}W_2(\varrho_s,\varrho'_s).
\]
Since \((\mathcal P_2(\mathbb R^n),W_2)\) is complete,
\((\mathcal M_t,d_t)\) is complete.

For \(\varrho\in\mathcal M_t\), the quantity
\(
  \sup_{s\in[t,T]}m_2(\varrho_s)
\)
is finite. Consider the frozen-law jump--diffusion system
\begin{align}
  \mathrm dX_s^\varrho
  &=
  f(X_s^\varrho,U_s^\varrho)\,\mathrm ds
  +
  \sigma(X_s^\varrho,U_s^\varrho)\,\mathrm dW_s,
  \label{eq:frozen-law-state-equation}
  \\
  U_s^\varrho
  &=
  V
  +
  \int_{(t,s]\times[0,1]}
    \bigl(
      \alpha_r(X_{r-}^\varrho,\varrho_r,z)-U_{r-}^\varrho
    \bigr)
    N(\mathrm dr,\mathrm dz),
  \qquad s\in[t,T],
  \label{eq:frozen-law-control-equation}
\end{align}
with \(X_t^\varrho=\xi\), where \((\xi,V)\) is square integrable with law
\(\zeta\) and independent of \((W,N)\).

The coefficient \( \mathfrak j_\varrho(s,x,u,z):=\alpha_s(x,\varrho_s,z)-u \)
is Borel measurable and satisfies
\[
  \int_0^1
    \|\mathfrak j_\varrho(s,x,u,z)-\mathfrak j_\varrho(s,x',u',z)\|^2\,\mathrm dz
  \le
  C\bigl(\|x-x'\|^2+\|u-u'\|^2\bigr),
\]
and
\[
  \int_0^1
    \|\mathfrak j_\varrho(s,x,u,z)\|^2\,\mathrm dz
  \le
  C\bigl(1+\|x\|^2+\|u\|^2+m_2(\varrho_s)\bigr).
\]
These estimates are the required \(L^2\)-Lipschitz and linear-growth estimates
for the Brownian and Poisson integral terms. Applying the standard
existence-and-uniqueness theorem for SDEs with Brownian and Poisson random
measure terms therefore gives a unique strong solution
\((X^\varrho,U^\varrho)\) \cite{Applebaum2009Levy,Situ2005SDEJumps}, and
\[
  \mathbb E\!\left[
    \sup_{s\in[t,T]}
    \bigl(
      \|X_s^\varrho\|^2+\|U_s^\varrho\|^2
    \bigr)
  \right]
  <\infty .
\]

Define \( \Phi(\varrho)_s:=\mathcal L(X_s^\varrho)\), for \(s\in[t,T]\).
For \(t\le r\le s\le T\),
\[
  W_2^2\bigl(\Phi(\varrho)_s,\Phi(\varrho)_r\bigr)
  \le
  \mathbb E\bigl[\|X_s^\varrho-X_r^\varrho\|^2\bigr].
\]
Using the state equation, the linear growth of \(f\) and \(\sigma\), and the
moment bound for \((X^\varrho,U^\varrho)\), one obtains
\[
  \mathbb E\bigl[\|X_s^\varrho-X_r^\varrho\|^2\bigr]
  \le
  C|s-r|.
\]
Thus \(\Phi(\varrho)\in\mathcal M_t\).

Let \(\varrho,\varrho'\in\mathcal M_t\), and construct the corresponding
frozen-law solutions with the same initial condition \((\xi,V)\), Brownian
motion \(W\), and Poisson random measure \(N\). Set \(\Delta_s^X:=X_s^\varrho-X_s^{\varrho'}\), and \(\Delta_s^U:=U_s^\varrho-U_s^{\varrho'}\).
Using standard \(L^2\) estimates for jump--diffusions, together with the
Lipschitz assumptions on \(f\), \(\sigma\), and \(\alpha\), one obtains
\[
  F(s)
  \le
  C\int_t^s
    \Bigl(
      F(\ell)+W_2(\varrho_\ell,\varrho'_\ell)^2
    \Bigr)\,\mathrm d\ell,
  \qquad s\in[t,T],
\]
where
\[
  F(s)
  :=
  \mathbb E\!\left[
    \sup_{r\in[t,s]}
    \bigl(
      \|\Delta_r^X\|^2+\|\Delta_r^U\|^2
    \bigr)
  \right].
\]
Since \(
  W_2(\varrho_\ell,\varrho'_\ell)\le d_t(\varrho,\varrho')\),
this gives
\[
  F(s)
  \le
  C\int_t^s F(\ell)\,\mathrm d\ell
  +
  C(s-t)d_t(\varrho,\varrho')^2.
\]
Solving this integral inequality gives
\[
  F(s)
  \le
  C(s-t)e^{C(s-t)}d_t(\varrho,\varrho')^2.
\]
Hence
\[
\begin{aligned}
  d_t\bigl(\Phi(\varrho),\Phi(\varrho')\bigr)^2
  &=
  \sup_{s\in[t,T]}
  W_2^2\bigl(
    \mathcal L(X_s^\varrho),
    \mathcal L(X_s^{\varrho'})
  \bigr)
  \\
  &\le
  \sup_{s\in[t,T]}
  \mathbb E\!\left[
    \|X_s^\varrho-X_s^{\varrho'}\|^2
  \right]
  \\
  &\le
  \sup_{s\in[t,T]}F(s)
  \\
  &\le
  C(T-t)e^{C(T-t)}
  d_t(\varrho,\varrho')^2.
\end{aligned}
\]
After enlarging the constant \(C\), this becomes
\[
  d_t\bigl(\Phi(\varrho),\Phi(\varrho')\bigr)^2
  \le
  C(T-t)d_t(\varrho,\varrho')^2.
\]
Therefore \(\Phi\) is a contraction whenever \(T-t\) is sufficiently small
\cite{Sznitman1991Propagation}. Consequently, on sufficiently short intervals,
Banach's fixed-point theorem gives a unique local closed-loop solution.

We next prove an a priori moment estimate for any such local closed-loop
solution. Write the control equation in compensated form. Let
\( \widetilde N(\mathrm dr,\mathrm dz) := N(\mathrm dr,\mathrm dz)-\lambda\,\mathrm dr\,\mathrm dz\).
Then
\begin{align*}
  U_s
  &=
  U_t
  +
  \int_{(t,s]\times[0,1]}
    \bigl(
      \alpha_r(X_{r-},\mu_r,z)-U_{r-}
    \bigr)
    \widetilde N(\mathrm dr,\mathrm dz)
  \\
  &\quad+
  \lambda\int_t^s\!\int_0^1
    \bigl(
      \alpha_r(X_r,\mu_r,z)-U_r
    \bigr)\,\mathrm dz\,\mathrm dr .
\end{align*}
Using standard maximal \(L^2\) estimates for stochastic integrals, the linear
growth of \(f\) and \(\sigma\), the second-moment bound on \(\alpha\), and
\[
  m_2(\mu_s)
  \le
  \mathbb E\!\left[
    \sup_{\ell\in[t,s]}\|X_\ell\|^2
  \right],
\]
one obtains
\[
  \mathbb E\!\left[
    \sup_{\ell\in[t,s]}
    \bigl(
      \|X_\ell\|^2+\|U_\ell\|^2
    \bigr)
  \right]
  \le
  C\Bigl(
    1+m_2(\zeta)
    +
    \int_t^s
      \mathbb E\!\left[
        \sup_{v\in[t,\ell]}
        \bigl(
          \|X_v\|^2+\|U_v\|^2
        \bigr)
      \right]\mathrm d\ell
  \Bigr).
\]
Solving this integral inequality gives
\[
  \mathbb E\!\left[
    \sup_{\ell\in[t,T]}
    \bigl(
      \|X_\ell\|^2+\|U_\ell\|^2
    \bigr)
  \right]
  \le
  C\bigl(1+m_2(\zeta)\bigr)
\]
on every interval on which the local solution is constructed.

We now extend the local construction to the whole interval \([t,T]\). Choose a
partition \( t=t_0<t_1<\cdots<t_N=T \) such that
\[
  C(t_{i+1}-t_i)<1,
  \qquad i=0,\ldots,N-1,
\]
where \(C\) is the constant in the local contraction estimate. This constant is
uniform: it depends only on the coefficient bounds, the refresh rate,
\(L_\alpha\), \(M_\alpha\), and the horizon \(T\), and not on the particular
initial law used on a subinterval. Hence such a finite partition exists. On
\([t_0,t_1]\), the contraction argument gives a unique local closed-loop
solution. Let \( \eta_{t_1}:=\mathcal L(X_{t_1},U_{t_1})\).
The moment estimate implies that \(m_2(\eta_{t_1})<\infty\). Therefore
\(\eta_{t_1}\) is an admissible initial joint law for the next subinterval.
Repeating the same fixed-point argument on \([t_1,t_2]\) with the restricted
sampling rule \( \alpha|_{[t_1,t_2]}\), and iterating over the partition,
constructs a weak closed-loop solution on the whole interval \([t,T]\). The
moment estimate prevents explosion of the restarted initial laws, so the
iteration reaches \(T\).

We now prove uniqueness of the \((X,U)\)-law. Let \((X,U)\) be any weak
closed-loop solution for \((t,\zeta,\alpha)\), and let
\( \mu_s:=\mathcal L(X_s) \), for \(s\in[t,T]\). Then \((X,U)\) solves the
frozen-law system
\eqref{eq:frozen-law-state-equation}--%
\eqref{eq:frozen-law-control-equation}
with \(\varrho=\mu\). Hence \(\mu=\Phi(\mu)\) on each subinterval in the above
partition. Since the fixed point is unique on each subinterval, the state-law
flow is unique. Once the state-law flow is fixed, the coefficients in the
closed-loop equations are deterministic functions of \((s,x,u,z)\). Uniqueness
in law for the resulting classical jump--diffusion then yields uniqueness of
the \((X,U)\)-law on \(\mathcal Y_{t,T}\). We denote this unique law by
\(\mathbb P^{t,\zeta,\alpha}\).

The preceding local moment estimate therefore extends to the global solution on
\([t,T]\), which gives \eqref{eq:moment-bound-joint-poisson}.

To prove stability in the initial joint law, let
\(\zeta,\zeta'\in\mathcal P_2(\mathbb R^n\times\mathbb R^q)\). Choose
square-integrable random vectors \((\xi,V)\), \((\xi',V')\) on a common
probability space such that
\[
  \mathcal L(\xi,V)=\zeta,
  \qquad
  \mathcal L(\xi',V')=\zeta',
  \qquad
  \mathbb E\!\left[
    \|\xi-\xi'\|^2+\|V-V'\|^2
  \right]
  =
  W_2(\zeta,\zeta')^2 .
\]
Let \((X,U)\) and \((X',U')\) be the closed-loop solutions driven by the same
\((W,N)\) and initialized by \((\xi,V)\) and \((\xi',V')\), respectively.
Denoting expectation on this space by \(\mathbb E\), set
\[
  G(s):=
  \mathbb E\!\left[
    \sup_{r\in[t,s]}
    \bigl(
      \|X_r-X_r'\|^2+\|U_r-U_r'\|^2
    \bigr)
  \right].
\]
Arguing as above,
\[
  G(s)
  \le
  C W_2(\zeta,\zeta')^2
  +
  C\int_t^s
    W_2\bigl(
      \mu_\ell^{t,\zeta,\alpha},
      \mu_\ell^{t,\zeta',\alpha}
    \bigr)^2\,\mathrm d\ell.
\]
Since the state-law distance is bounded by the corresponding joint-law
distance, the previous estimate gives
\[
  G(s)
  \le
  C W_2(\zeta,\zeta')^2
  +
  C\int_t^s G(\ell)\,\mathrm d\ell.
\]
Solving this integral inequality gives
\[
  G(s)\le C W_2(\zeta,\zeta')^2,
  \qquad s\in[t,T].
\]
Since
\[
  W_2^2\bigl(
    \eta_\ell^{t,\zeta,\alpha},
    \eta_\ell^{t,\zeta',\alpha}
  \bigr)
  \le
  \mathbb E\!\left[
    \|X_\ell-X_\ell'\|^2+\|U_\ell-U_\ell'\|^2
  \right],
\]
we obtain
\[
  \sup_{\ell\in[t,T]}
  W_2\bigl(
    \eta_\ell^{t,\zeta,\alpha},
    \eta_\ell^{t,\zeta',\alpha}
  \bigr)
  \le
  C W_2(\zeta,\zeta').
\]

For \(t\le r\le s\le T\),
\[
  W_2^2\bigl(
    \mu_s^{t,\zeta,\alpha},
    \mu_r^{t,\zeta,\alpha}
  \bigr)
  \le
  \mathbb E^{t,\zeta,\alpha}\bigl[
    \|X_s-X_r\|^2
  \bigr].
\]
Using the state equation, the linear growth of \(f\) and \(\sigma\), and the
moment bound already proved, one obtains
\[
  \mathbb E^{t,\zeta,\alpha}\bigl[
    \|X_s-X_r\|^2
  \bigr]
  \le
  C\bigl(1+m_2(\zeta)\bigr)|s-r|.
\]
Therefore
\[
  W_2\bigl(
    \mu_s^{t,\zeta,\alpha},
    \mu_r^{t,\zeta,\alpha}
  \bigr)
  \le
  C\bigl(1+m_2(\zeta)^{1/2}\bigr)|s-r|^{1/2},
\]
and \(s\mapsto\mu_s^{t,\zeta,\alpha}\) is continuous in \(W_2\).

Finally, we prove the continuity of
\(s\mapsto\eta_s^{t,\zeta,\alpha}\). If \(s_k\downarrow s\), then, by the
c\`adl\`ag paths of \(U\) and the continuity of \(X\),
\((X_{s_k},U_{s_k})\to(X_s,U_s)\) a.s. Together with the uniform second-moment
bound, this implies
\[
  W_2\bigl(
    \eta_{s_k}^{t,\zeta,\alpha},
    \eta_s^{t,\zeta,\alpha}
  \bigr)\to0.
\]
If \(s_k\uparrow s\), then
\((X_{s_k},U_{s_k})\to(X_s,U_{s-})\) a.s. Hence
\[
  \eta_{s_k}^{t,\zeta,\alpha}
  \longrightarrow
  \mathcal L(X_s,U_{s-})
  \quad\text{in }W_2.
\]
Since the Poisson random measure has no atom at any fixed deterministic time
\(s\), we have \(U_{s-}=U_s\) almost surely. Therefore
\[
  \mathcal L(X_s,U_{s-})=\mathcal L(X_s,U_s)
  =
  \eta_s^{t,\zeta,\alpha}.
\]
Thus
\[
  s\longmapsto\eta_s^{t,\zeta,\alpha}
  \in
  C\bigl([t,T];\mathcal P_2(\mathbb R^n\times\mathbb R^q)\bigr).
\]
\end{proof}


For the initial problem at time \(0\), the natural initial datum is a joint law
\(\zeta_0\in\mathcal P_2(\mathbb R^n\times\mathbb R^q)\). If one wishes to start only from a prescribed initial state law
\(\mu_0\in\mathcal P_2(\mathbb R^n)\), then an initial control law must also be
specified. Equivalently, one chooses a Borel kernel
\(\rho_0(\mathrm du\mid x,\mu_0)\) on \(\mathbb R^q\) satisfying
\(
  \int_{\mathbb R^n}\int_{\mathbb R^q}
    \|u\|^2\,\rho_0(\mathrm du\mid x,\mu_0)\,\mu_0(\mathrm dx)
  <\infty,
\)
and sets
\[
  \zeta_0(\mathrm dx,\mathrm du)
  :=
  \mu_0(\mathrm dx)\rho_0(\mathrm du\mid x,\mu_0).
\]
In that case we write
\[
  \mathbb P^\alpha:=\mathbb P^{0,\zeta_0,\alpha},
  \qquad
  \eta_s^\alpha:=\eta_s^{0,\zeta_0,\alpha},
  \qquad
  \mu_s^\alpha:=\mu_s^{0,\zeta_0,\alpha},
  \qquad
  \alpha\in\mathfrak A.
\]

We next record a structural property that will be used in the dynamic
programming argument in Section \ref{sec:dpp-poisson}, namely that admissible sampling rules may be restricted to subintervals and concatenated in time.

\begin{proposition}\label{prop:stability-sampling-rules}
Let \(0\le t\le \tau\le T\).
\begin{enumerate}
\item
If \(\alpha\in\mathfrak A[t,T]\), then
\(
  \alpha|_{[t,\tau]}\in\mathfrak A[t,\tau]
\)
and
\(
  \alpha|_{[\tau,T]}\in\mathfrak A[\tau,T].
\)

\item
If
\(
  \alpha^1\in\mathfrak A[t,\tau]
\)
and
\(
  \alpha^2\in\mathfrak A[\tau,T],
\)
then the concatenated sampling rule
\[
  (\alpha^1\otimes_\tau \alpha^2)_s(x,\mu,z)
  :=
  \begin{cases}
    \alpha_s^1(x,\mu,z), & s\in[t,\tau],\\
    \alpha_s^2(x,\mu,z), & s\in(\tau,T],
  \end{cases}
\]
belongs to \(\mathfrak A[t,T]\).
\end{enumerate}
\end{proposition}

\begin{proof}
Restriction plainly preserves Borel measurability and the estimates
\eqref{eq:sampling-rule-lipschitz}--\eqref{eq:sampling-rule-second-moment},
with the same constants \(L_\alpha\) and \(M_\alpha\).

For concatenation, the map
\((s,x,\mu,z)\mapsto(\alpha^1\otimes_\tau\alpha^2)_s(x,\mu,z)\) is Borel
measurable because it is defined by two Borel maps on the Borel sets
\([t,\tau]\) and \((\tau,T]\). On \([t,\tau]\), the estimates hold because
\(\alpha^1\in\mathfrak A[t,\tau]\); on \((\tau,T]\), they hold because
\(\alpha^2\in\mathfrak A[\tau,T]\). Since both classes are defined using the
same fixed constants \(L_\alpha\) and \(M_\alpha\), the concatenated rule
satisfies the same estimates on the whole interval \([t,T]\). Hence
\(\alpha^1\otimes_\tau\alpha^2\in\mathfrak A[t,T]\).
\end{proof}

The next proposition records the corresponding restriction and concatenation
properties of the closed-loop laws. Since the future dynamics depend on both
the current state \(X\) and the current control value \(U\), the natural
continuation variable is the joint law of \((X,U)\), not only the law of \(X\) \cite{DjetePossamaiTan2022MKVDPP}.

\begin{proposition}\label{prop:stability-joint-laws}
Let \(0\le t\le\tau\le T\).
\begin{enumerate}
\item
If \(\alpha\in\mathfrak A[t,T]\), then
\begin{equation}
  (\Pi_{t,\tau})_\#\mathbb P^{t,\zeta,\alpha}
  =
  \mathbb P^{t,\zeta,\alpha|_{[t,\tau]}}.
  \label{eq:restriction-path-law-poisson}
\end{equation}

\item
Let \(\alpha^1\in\mathfrak A[t,\tau]\) and
\(\alpha^2\in\mathfrak A[\tau,T]\), and set
\(
  \vartheta:=\eta_\tau^{t,\zeta,\alpha^1}\in
  \mathcal P_2(\mathbb R^n\times\mathbb R^q).
\)
Then the joint law flow of the concatenated rule
\(\alpha^1\otimes_\tau\alpha^2\) is given by
\begin{equation}
  \eta_s^{t,\zeta,\alpha^1\otimes_\tau\alpha^2}
  =
  \begin{cases}
    \eta_s^{t,\zeta,\alpha^1}, & s\in[t,\tau],\\
    \eta_s^{\tau,\vartheta,\alpha^2}, & s\in[\tau,T],
  \end{cases}
  \label{eq:concatenated-joint-law-flow}
\end{equation}
and therefore
\begin{equation}
  \mu_s^{t,\zeta,\alpha^1\otimes_\tau\alpha^2}
  =
  \begin{cases}
    \mu_s^{t,\zeta,\alpha^1}, & s\in[t,\tau],\\
    \mu_s^{\tau,\vartheta,\alpha^2}, & s\in[\tau,T].
  \end{cases}
  \label{eq:concatenated-state-law-flow}
\end{equation}
Consequently,
\begin{align}
  (\Pi_{t,\tau})_\#\mathbb P^{t,\zeta,\alpha^1\otimes_\tau\alpha^2}
  &=
  \mathbb P^{t,\zeta,\alpha^1},
  \label{eq:concat-first-restriction-poisson}
  \\
  (\Pi_{\tau,T})_\#\mathbb P^{t,\zeta,\alpha^1\otimes_\tau\alpha^2}
  &=
  \mathbb P^{\tau,\vartheta,\alpha^2}.
  \label{eq:concat-second-restriction-poisson}
\end{align}
\end{enumerate}
\end{proposition}

\begin{proof}
For part~(1), let \((X,U,W,N)\) be a weak closed-loop solution for
\((t,\zeta,\alpha)\) on some filtered probability space. Restrict the processes
to \([t,\tau]\). Since \(\alpha|_{[t,\tau]}\) agrees with \(\alpha\) on
\([t,\tau]\), the restricted quadruple is a weak closed-loop solution for
\((t,\zeta,\alpha|_{[t,\tau]})\). By uniqueness in law from
Theorem~\ref{thm:wellposedness-poisson-refresh}, the law of the restricted
process is exactly \(\mathbb P^{t,\zeta,\alpha|_{[t,\tau]}}\), which proves
\eqref{eq:restriction-path-law-poisson}.

For part~(2), let \((X,U,W,N)\) be a weak closed-loop solution for
\((t,\zeta,\alpha^1\otimes_\tau\alpha^2)\), and let
\(
  \mathbb P:=\mathbb P^{t,\zeta,\alpha^1\otimes_\tau\alpha^2}.
\)
By part~(1),
\[
  (\Pi_{t,\tau})_\#\mathbb P
  =
  \mathbb P^{t,\zeta,\alpha^1},
\]
hence
\[
  \eta_s^{t,\zeta,\alpha^1\otimes_\tau\alpha^2}
  =
  \eta_s^{t,\zeta,\alpha^1},
  \qquad
  \mu_s^{t,\zeta,\alpha^1\otimes_\tau\alpha^2}
  =
  \mu_s^{t,\zeta,\alpha^1},
  \qquad
  s\in[t,\tau],
\]
and in particular
\[
  \eta_\tau^{t,\zeta,\alpha^1\otimes_\tau\alpha^2}
  =
  \eta_\tau^{t,\zeta,\alpha^1}
  =:\vartheta.
\]

Now consider the shifted processes on \([\tau,T]\):
\[
  \bar W_s:=W_s-W_\tau,
  \qquad
  \bar N((\tau,s]\times B):=N((\tau,s]\times B),
  \qquad
  s\in[\tau,T].
\]
Since \((X_\tau,U_\tau)\) is \(\mathcal F_\tau\)-measurable and the increments
of \(W\) and \(N\) after time \(\tau\) are independent of \(\mathcal F_\tau\),
\((X_\tau,U_\tau)\) is independent of \((\bar W,\bar N)\). In addition, since
\(\alpha^1\otimes_\tau\alpha^2\) agrees with \(\alpha^2\) on \((\tau,T]\), the
restricted quadruple
\[
  \bigl(
    X|_{[\tau,T]},\,
    U|_{[\tau,T]},\,
    \bar W,\,
    \bar N
  \bigr)
\]
is a weak closed-loop solution for \((\tau,\vartheta,\alpha^2)\). By uniqueness
in law from Theorem~\ref{thm:wellposedness-poisson-refresh},
\(
  (\Pi_{\tau,T})_\#\mathbb P
  =
  \mathbb P^{\tau,\vartheta,\alpha^2}.
\)
This yields \eqref{eq:concatenated-joint-law-flow},
\eqref{eq:concatenated-state-law-flow}, and hence
\eqref{eq:concat-first-restriction-poisson}--%
\eqref{eq:concat-second-restriction-poisson}.
\end{proof}

The optimization problem studied in the sequel will be defined after the local
control and state costs have been introduced. In the present Poisson-refresh
formulation, the fundamental closed-loop object is the path law
\(\mathbb P^{t,\zeta,\alpha}\) of the joint process \((X,U)\), with associated
joint law flow \(\eta^{t,\zeta,\alpha}\) and state-law flow
\(\mu^{t,\zeta,\alpha}\).

\section{Reference sampling rule and local control entropy}
\label{sec:control-entropy}
We now introduce the first local cost in the Poisson-refresh formulation of
Section~\ref{sec:poisson-refresh-formulation}, namely the control-side entropy cost.
Because the control is represented by an actual sampled process \(U\), this
cost is naturally associated with the refresh mechanism itself. It quantifies,
locally in time and state, how the sampling rule differs from a prescribed
reference sampling rule.

To express this local comparison, we pass from the sampled process itself to
the distribution generated at a refresh time. The relevant object is therefore
the instantaneous law of the post-refresh control value, conditional on the
current state \(x\) and the current state law \(\mu\). This focuses the cost on
the control update generated at refresh times, rather than on the full joint
path law of \((X,U)\), whose relative entropy would also reflect the induced
evolution of the state process. As in
Section~\ref{sec:poisson-refresh-formulation}, this conditional law is
determined by the sampling rule through its induced feedback kernel.

\medskip

Consequently, we fix throughout a reference sampling rule
\(
  \alpha^I=(\alpha_t^I)_{t\in[0,T]}
  \in
  \mathfrak A,
\)
and denote by
\[
  c_t^I(\cdot\mid x,\mu)
  :=
  c_t^{\alpha^I}(\cdot\mid x,\mu),
  \qquad
  (t,x,\mu)\in
  [0,T]\times\mathbb R^n\times\mathcal P_2(\mathbb R^n),
\]
the induced reference feedback kernel from
Proposition~\ref{prop:kernel-induced-by-sampling-rule}, that is,
\begin{equation}
  c_t^I(B\mid x,\mu)
  :=
  \int_0^1
    \mathbf 1_B\!\bigl(
      \alpha_t^I(x,\mu,z)
    \bigr)\,\mathrm dz,
  \qquad
  B\in\mathcal B(\mathbb R^q).
  \label{eq:reference-kernel-from-reference-sampling-rule}
\end{equation}
By Proposition~\ref{prop:kernel-induced-by-sampling-rule}, the family
\(c^I=(c_t^I)_{t\in[0,T]}\) is Borel measurable and satisfies the same
Lipschitz and second-moment bounds as any kernel induced by an admissible
sampling rule.

Let \(\mathcal K\) denote the class of all Borel probability kernels
\(\gamma:\mathbb R^n\to\mathcal P(\mathbb R^q)\), written
\(x\mapsto\gamma(\cdot\mid x)\).
In the local comparison below, the state law \(\mu\) is fixed. Thus
\(\gamma\) is viewed as a kernel from the current state \(x\) to a candidate
post-refresh control law. When such a kernel is induced by an admissible
sampling rule, it is obtained by evaluating
\(c^\alpha_t(\cdot\mid x,\mu)\) at this fixed \(\mu\).

For \(t\in[0,T]\), \(\mu\in\mathcal P_2(\mathbb R^n)\), and
\(\gamma\in\mathcal K\), define
\begin{align}
  M_{t,\mu}^\gamma(\mathrm dx,\mathrm du)
  &:=
  \mu(\mathrm dx)\,\gamma(\mathrm du\mid x),
  \\
  M_{t,\mu}^I(\mathrm dx,\mathrm du)
  &:=
  \mu(\mathrm dx)\,c_t^I(\mathrm du\mid x,\mu).
\end{align}
Here the second coordinate \(u\) is interpreted as a candidate
post-refresh value, not as the current value of the control process.
The Poisson refresh rate \(\lambda>0\) was fixed in
Assumption~\ref{assump:primitive-coeff}. For
\((t,x,\mu,\rho)\in
[0,T]\times\mathbb R^n\times\mathcal P_2(\mathbb R^n)\times\mathcal P(\mathbb R^q)\),
define the pointwise control-cost rate
\begin{equation}
  \kappa_{\mathrm{ctrl}}(t,x,\mu,\rho)
  :=
  \lambda\,
  D\!\Bigl(
    \rho\,\Big\|\,c_t^I(\cdot\mid x,\mu)
  \Bigr),
  \label{eq:control-KL-rate-poisson}
\end{equation}
with the convention that the value is \(+\infty\) whenever
\(
  \rho\not\ll c_t^I(\cdot\mid x,\mu)
\).

Since relative entropy is lower semicontinuous, hence Borel measurable, on
\(
  \mathcal P(\mathbb R^q)\times\mathcal P(\mathbb R^q)
\)
for the topology of weak convergence, and since
\(
  (t,x,\mu)\longmapsto c_t^I(\cdot\mid x,\mu)
\)
is Borel measurable, it follows that \(\kappa_{\mathrm{ctrl}}\) is Borel
measurable on
\(
  [0,T]\times\mathbb R^n\times\mathcal P_2(\mathbb R^n)\times\mathcal P(\mathbb R^q).
\)
Moreover, for each fixed \((t,x,\mu)\), the map
\[
  \rho\longmapsto \kappa_{\mathrm{ctrl}}(t,x,\mu,\rho)
\]
is convex and lower semicontinuous on \(\mathcal P(\mathbb R^q)\).

For \(t\in[0,T]\), \(\mu\in\mathcal P_2(\mathbb R^n)\), and
\(\gamma\in\mathcal K\), define the local control functional
\begin{equation}
  \mathcal J_{\mathrm{ctrl}}(t,\mu;\gamma)
  :=
  \int_{\mathbb R^n}
    \kappa_{\mathrm{ctrl}}\!\bigl(
      t,x,\mu,\gamma(\cdot\mid x)
    \bigr)\,
    \mu(\mathrm dx)
  \in [0,\infty].
  \label{eq:local-control-functional-poisson}
\end{equation}

The next lemma identifies \(\mathcal J_{\mathrm{ctrl}}\) with a relative entropy
on the joint state--refresh-value space.

\begin{lemma}\label{lem:joint-entropy-disintegration-poisson}
For every \(t\in[0,T]\), every \(\mu\in\mathcal P_2(\mathbb R^n)\), and every
\(\gamma\in\mathcal K\),
\begin{equation}
  \mathcal J_{\mathrm{ctrl}}(t,\mu;\gamma)
  =
  \lambda\,
  D\!\Bigl(
    M_{t,\mu}^\gamma\,\Big\|\,M_{t,\mu}^I
  \Bigr).
  \label{eq:joint-entropy-identity-poisson}
\end{equation}
Consequently, for each fixed \((t,\mu)\), the map
\(
  \gamma\longmapsto \mathcal J_{\mathrm{ctrl}}(t,\mu;\gamma)
\)
is convex, and
\[
  \mathcal J_{\mathrm{ctrl}}(t,\mu;\gamma)=0
  \quad\Longleftrightarrow\quad
  \gamma(\cdot\mid x)=c_t^I(\cdot\mid x,\mu)
  \text{ for }\mu\text{-a.e.\ }x.
\]
In particular,
\begin{equation}
  \mathcal J_{\mathrm{ctrl}}\!\bigl(
    t,\mu;c_t^I(\cdot\mid\cdot,\mu)
  \bigr)
  =
  0.
  \label{eq:local-control-zero-reference-poisson}
\end{equation}
\end{lemma}

\begin{proof}
The measures \(M_{t,\mu}^\gamma\) and \(M_{t,\mu}^I\) have the same first
marginal \(\mu\). The chain rule for relative entropy with respect to this
common disintegration therefore yields \cite{DemboZeitouni1998LDP}
\[
  D\!\Bigl(
    M_{t,\mu}^\gamma\,\Big\|\,M_{t,\mu}^I
  \Bigr)
  =
  \int_{\mathbb R^n}
    D\!\Bigl(
      \gamma(\cdot\mid x)\,\Big\|\,c_t^I(\cdot\mid x,\mu)
    \Bigr)\,
    \mu(\mathrm dx),
\]
with values in \([0,\infty]\). Multiplying by \(\lambda\) gives
\eqref{eq:joint-entropy-identity-poisson}, and the remaining assertions follow
immediately from the corresponding basic properties of relative entropy.
\end{proof}

The local control cost is defined at the level of the induced kernel rather than
at the level of a particular sampling parametrization. Indeed, different
sampling rules may push forward Lebesgue measure on \([0,1]\) to the same
post-refresh control law. The control-side cost therefore depends only on this
law, namely on the kernel \(\gamma\), and not on the particular choice of
auxiliary parametrization by \(z\).

We now evaluate this local functional along the closed-loop law generated by an
admissible sampling rule. Fix \(0\le t\le T\),
\(\zeta\in\mathcal P_2(\mathbb R^n\times\mathbb R^q)\), and
\(\alpha\in\mathfrak A[t,T]\). Let
\(
  \mu_s^{t,\zeta,\alpha}
  :=
  \mathcal L_{\mathbb P^{t,\zeta,\alpha}}(X_s),\) for every \(
  s\in[t,T],\)
be the state-law flow from
Theorem~\ref{thm:wellposedness-poisson-refresh}. For each \(s\in[t,T]\), define
the instantaneous closed-loop refresh kernel
\begin{equation}
  \gamma_s^{t,\zeta,\alpha}(\cdot\mid x)
  :=
  c_s^\alpha(\cdot\mid x,\mu_s^{t,\zeta,\alpha}),
  \qquad
  x\in\mathbb R^n.
  \label{eq:gamma-from-alpha-poisson}
\end{equation}
Then \(x\mapsto\gamma_s^{t,\zeta,\alpha}(\cdot\mid x)\) belongs to
\(\mathcal K\) for every \(s\in[t,T]\).

The control contribution of \(\alpha\) on \([t,T]\) is defined by
\begin{align}
  \mathcal H_{\mathrm{ctrl}}(t,\zeta;\alpha)
  &:=
  \int_t^T
    \mathcal J_{\mathrm{ctrl}}\!\bigl(
      s,\mu_s^{t,\zeta,\alpha};\gamma_s^{t,\zeta,\alpha}
    \bigr)\,\mathrm ds
  \label{eq:control-cost-general-time-poisson}
  \\
  &=
  \mathbb E^{t,\zeta,\alpha}\!\left[
    \int_t^T
      \kappa_{\mathrm{ctrl}}\!\bigl(
        s,X_s,\mu_s^{t,\zeta,\alpha},
        c_s^\alpha(\cdot\mid X_s,\mu_s^{t,\zeta,\alpha})
      \bigr)\,\mathrm ds
  \right],
  \label{eq:control-KL-rate-integrated-poisson}
\end{align}
with values in \([0,\infty]\). The equality of
\eqref{eq:control-cost-general-time-poisson} and
\eqref{eq:control-KL-rate-integrated-poisson} follows by expanding the
definition of \(\mathcal J_{\mathrm{ctrl}}\), using
\(\mu_s^{t,\zeta,\alpha}=\mathcal L_{\mathbb P^{t,\zeta,\alpha}}(X_s)\), and
interchanging the time and state integrations. This functional measures only the entropy cost of the refresh mechanism. It is
therefore not the relative entropy between the full path law
\(\mathbb P^{t,\zeta,\alpha}\) of \((X,U)\) and the corresponding path law under
the reference sampling rule, since such a path-space entropy would also include
the induced change in the state dynamics.

Since the reference sampling
rule \(\alpha^I\) induces the reference kernel \(c^I\), one has
\(
  \mathcal H_{\mathrm{ctrl}}\!\bigl(
    t,\zeta;\alpha^I|_{[t,T]}
  \bigr)
  =0
\)
for every \((t,\zeta)\).

The next proposition records an auxiliary marked-refresh representation of the
control-side entropy. Because refreshes may redraw the current control value,
the latent refresh times are in general not recoverable from the observed
control path. The representation is therefore formulated on an augmented
marked-refresh space, rather than on \(D([t,T];\mathbb R^q)\) alone. 

For \(0\le t\le T\), let \(\mathfrak M_{t,T}^q\) denote the space of finite
point measures on \((t,T]\times\mathbb R^q\), equipped with its standard Borel
\(\sigma\)-field. Set
\[
  \widehat{\mathcal U}_{t,T}
  :=
  \mathbb R^q\times \mathfrak M_{t,T}^q.
\]
A point \((u, \mathfrak m)\in\widehat{\mathcal U}_{t,T}\) is interpreted as an initial
control value \(u\) together with a marked refresh point measure
\[
   \mathfrak m=\sum_{k=1}^N\delta_{(\tau_k,v_k)},
\]
where \(N\in\mathbb N_0\), \(t<\tau_1<\cdots<\tau_N\le T\), and
\(v_k\in\mathbb R^q\) are the successive post-refresh values.

\begin{proposition}[Auxiliary marked-refresh representation of the control entropy]
\label{prop:control-KL-poisson-augmented}
Fix \(0\le t\le T\), \(\zeta\in\mathcal P_2(\mathbb R^n\times\mathbb R^q)\),
and \(\alpha\in\mathfrak A[t,T]\). For each
\((\mathbf x,u)\in\mathcal X_{t,T}\times\mathbb R^q\), define two auxiliary
marked-refresh laws
\[
  \widehat{\mathbf Q}^{\alpha;\mathbf x,u},
  \qquad
  \widehat{\mathbf Q}^{I;\mathbf x,u}
  \in
  \mathcal P(\widehat{\mathcal U}_{t,T})
\]
as follows:
\begin{enumerate}
\item
the initial control value is fixed and equal to \(u\) under both laws;

\item
the refresh times are given by a Poisson process of rate \(\lambda\), common to
both laws;

\item
conditional on the ordered refresh times
\(
  t<\tau_1<\tau_2<\cdots<\tau_N\le T,
\)
the post-refresh values are conditionally independent, and the conditional law
of the \(k\)-th post-refresh value under
\(\widehat{\mathbf Q}^{\alpha;\mathbf x,u}\) is
\(
  c^\alpha_{\tau_k}\!\bigl(
    \cdot\mid \mathbf x(\tau_k),\mu_{\tau_k}^{t,\zeta,\alpha}
  \bigr),
\)
whereas under \(\widehat{\mathbf Q}^{I;\mathbf x,u}\) it is
\(
  c^I_{\tau_k}\!\bigl(
    \cdot\mid \mathbf x(\tau_k),\mu_{\tau_k}^{t,\zeta,\alpha}
  \bigr).
\)
\end{enumerate}
These are auxiliary laws defined from the frozen path \(\mathbf x\); no claim is
made here that they are regular conditional laws of the actual latent refresh
mechanism given the realised state path.

The maps
\[
  (\mathbf x,u)\longmapsto \widehat{\mathbf Q}^{\alpha;\mathbf x,u},
  \qquad
  (\mathbf x,u)\longmapsto \widehat{\mathbf Q}^{I;\mathbf x,u}
\]
may be chosen Borel measurable from
\(\mathcal X_{t,T}\times\mathbb R^q\) into
\(\mathcal P(\widehat{\mathcal U}_{t,T})\). Moreover, for every
\((\mathbf x,u)\in\mathcal X_{t,T}\times\mathbb R^q\),
\begin{equation}
  D\!\Bigl(
    \widehat{\mathbf Q}^{\alpha;\mathbf x,u}
    \,\Big\|\,
    \widehat{\mathbf Q}^{I;\mathbf x,u}
  \Bigr)
  =
  \lambda
  \int_t^T
    D\!\Bigl(
      c_s^\alpha(\cdot\mid \mathbf x(s),\mu_s^{t,\zeta,\alpha})
      \,\Big\|\,
      c_s^I(\cdot\mid \mathbf x(s),\mu_s^{t,\zeta,\alpha})
    \Bigr)\,
    \mathrm ds,
  \label{eq:control-poisson-pathwise-augmented}
\end{equation}
with values in \([0,\infty]\). Let
\(
  \Lambda^{t,\zeta,\alpha}
  :=
  \mathcal L_{\mathbb P^{t,\zeta,\alpha}}(X,U_t)
  \in \mathcal P(\mathcal X_{t,T}\times\mathbb R^q).
\)
Then
\begin{equation}
  \mathcal H_{\mathrm{ctrl}}(t,\zeta;\alpha)
  =
  \int_{\mathcal X_{t,T}\times\mathbb R^q}
    D\!\Bigl(
      \widehat{\mathbf Q}^{\alpha;\mathbf x,u}
      \,\Big\|\,
      \widehat{\mathbf Q}^{I;\mathbf x,u}
    \Bigr)\,
    \Lambda^{t,\zeta,\alpha}(\mathrm d\mathbf x,\mathrm du).
  \label{eq:control-path-entropy-augmented}
\end{equation}
\end{proposition}

\begin{proof}
For fixed \((t,\zeta,\alpha)\), the maps
\(
  (s,\mathbf x)\longmapsto
  c_s^\alpha(\cdot\mid \mathbf x(s),\mu_s^{t,\zeta,\alpha})\), and
\(  (s,\mathbf x)\longmapsto
  c_s^I(\cdot\mid \mathbf x(s),\mu_s^{t,\zeta,\alpha})
\)
are Borel kernels from \([t,T]\times\mathcal X_{t,T}\) into
\(\mathcal P(\mathbb R^q)\). This follows from the Borel measurability of the
evaluation map \((s,\mathbf x)\mapsto \mathbf x(s)\), the measurability of the
refresh kernels, and the measurability of
\(s\mapsto\mu_s^{t,\zeta,\alpha}\) \cite{Kallenberg2021FMP}.

For each \((\mathbf x,u)\), construct the marked point process by first drawing
a Poisson number of refresh times on \((t,T]\) with rate \(\lambda\), then,
conditional on the ordered times \((\tau_1,\ldots,\tau_N)\), drawing the marks
independently from the corresponding kernels. This construction gives, for
every cylinder set in \(\widehat{\mathcal U}_{t,T}\), a probability depending
Borel measurably on \((\mathbf x,u)\). Since the finite point-measure space is
equipped with its standard Borel \(\sigma\)-field, the monotone-class theorem
then gives Borel measurability of
\[
  (\mathbf x,u)\mapsto \widehat{\mathbf Q}^{\alpha;\mathbf x,u},
  \qquad
  (\mathbf x,u)\mapsto \widehat{\mathbf Q}^{I;\mathbf x,u}.
\]

Fix \((\mathbf x,u)\in\mathcal X_{t,T}\times\mathbb R^q\). Under
\(\widehat{\mathbf Q}^{\alpha;\mathbf x,u}\) and
\(\widehat{\mathbf Q}^{I;\mathbf x,u}\), the initial control value and the
refresh-time law are the same; only the conditional laws of the post-refresh
marks differ. Let \(N_{t,T}\) denote the number of refresh times in \((t,T]\).
Conditional on \(\{N_{t,T}=n\}\) and on the ordered refresh times
\[
  (\tau_1,\ldots,\tau_n)=(s_1,\ldots,s_n),
  \qquad
  t<s_1<\cdots<s_n\le T,
\]
the conditional law of the refresh-mark vector \((V_1,\ldots,V_n)\) under
\(\widehat{\mathbf Q}^{\alpha;\mathbf x,u}\) is
\[
  \bigotimes_{k=1}^n
  c_{s_k}^\alpha\!\bigl(
    \cdot\mid \mathbf x(s_k),\mu_{s_k}^{t,\zeta,\alpha}
  \bigr),
\]
whereas under \(\widehat{\mathbf Q}^{I;\mathbf x,u}\) it is
\[
  \bigotimes_{k=1}^n
  c_{s_k}^I\!\bigl(
    \cdot\mid \mathbf x(s_k),\mu_{s_k}^{t,\zeta,\alpha}
  \bigr).
\]
Since the initial value and the refresh-time law coincide under the two
measures, the relative entropy reduces to the conditional relative entropy of
the marks. Conditional on
\(
  N_{t,T}=n
\)
and
\(
  (\tau_1,\ldots,\tau_n)=(s_1,\ldots,s_n),
\)
the conditional Radon--Nikodym derivative of the mark law under
\(
  \widehat{\mathbf Q}^{\alpha;\mathbf x,u}
\)
with respect to the mark law under
\(
  \widehat{\mathbf Q}^{I;\mathbf x,u}
\)
is
\[
  \prod_{k=1}^n
  \frac{
    \mathrm d c_{s_k}^{\alpha}(
      \cdot\mid \mathbf x(s_k),\mu_{s_k}^{t,\zeta,\alpha}
    )
  }{
    \mathrm d c_{s_k}^{I}(
      \cdot\mid \mathbf x(s_k),\mu_{s_k}^{t,\zeta,\alpha}
    )
  }(v_k),
\]
whenever the first product measure is absolutely continuous with respect to the
second. Taking logarithms gives
\[
  \sum_{k=1}^n
  \log
  \frac{
    \mathrm d c_{s_k}^{\alpha}(
      \cdot\mid \mathbf x(s_k),\mu_{s_k}^{t,\zeta,\alpha}
    )
  }{
    \mathrm d c_{s_k}^{I}(
      \cdot\mid \mathbf x(s_k),\mu_{s_k}^{t,\zeta,\alpha}
    )
  }(v_k).
\]
Taking expectation under the conditional mark law therefore yields
\[
  \sum_{k=1}^n
  D\!\Bigl(
    c_{s_k}^{\alpha}(
      \cdot\mid \mathbf x(s_k),\mu_{s_k}^{t,\zeta,\alpha}
    )
    \,\Big\|\,
    c_{s_k}^{I}(
      \cdot\mid \mathbf x(s_k),\mu_{s_k}^{t,\zeta,\alpha}
    )
  \Bigr).
\]
Consequently,
\[
  D\!\Bigl(
    \widehat{\mathbf Q}^{\alpha;\mathbf x,u}
    \,\Big\|\,
    \widehat{\mathbf Q}^{I;\mathbf x,u}
  \Bigr)
  =
  \mathbb E\!\left[
    \sum_{\tau_k\in(t,T]}
      D\!\Bigl(
        c_{\tau_k}^{\alpha}(
          \cdot\mid \mathbf x(\tau_k),\mu_{\tau_k}^{t,\zeta,\alpha}
        )
        \,\Big\|\,
        c_{\tau_k}^{I}(
          \cdot\mid \mathbf x(\tau_k),\mu_{\tau_k}^{t,\zeta,\alpha}
        )
      \Bigr)
  \right].
\]
The expectation on the right-hand side is taken over the common auxiliary
Poisson refresh clock. Therefore, by the intensity formula for a Poisson process of
rate \(\lambda\),
\[
  \mathbb E\!\left[
    \sum_{\tau_k\in(t,T]}
      h(\tau_k)
  \right]
  =
  \lambda\int_t^T h(s)\,\mathrm ds
\]
for nonnegative measurable \(h\). Applying this with
\(
  h(s)
  =
  D\!\Bigl(
    c_s^\alpha(\cdot\mid \mathbf x(s),\mu_s^{t,\zeta,\alpha})
    \,\Big\|\,
    c_s^I(\cdot\mid \mathbf x(s),\mu_s^{t,\zeta,\alpha})
  \Bigr)
\)
proves \eqref{eq:control-poisson-pathwise-augmented}.

Finally, integrate \eqref{eq:control-poisson-pathwise-augmented} with respect
to the joint law of the state path \(X\) and the initial control value \(U_t\)
under \(\mathbb P^{t,\zeta,\alpha}\). Since the
right-hand side of \eqref{eq:control-poisson-pathwise-augmented} does not
depend on \(u\), this gives
\[
  \int_t^T
    \mathbb E^{t,\zeta,\alpha}\!\left[
      \lambda
      D\!\Bigl(
        c_s^\alpha(\cdot\mid X_s,\mu_s^{t,\zeta,\alpha})
        \,\Big\|\,
        c_s^I(\cdot\mid X_s,\mu_s^{t,\zeta,\alpha})
      \Bigr)
    \right]\mathrm ds.
\]
By the definition of \(\kappa_{\mathrm{ctrl}}\) and
\(\mathcal H_{\mathrm{ctrl}}\), this is exactly
\(\mathcal H_{\mathrm{ctrl}}(t,\zeta;\alpha)\). Hence
\eqref{eq:control-path-entropy-augmented} follows.
\end{proof}

Proposition \ref{prop:control-KL-poisson-augmented} represents the control entropy on an augmented marked-refresh space, where the refresh times and post-refresh values are kept
explicit. The actual closed-loop model, however, observes only the resulting
c\`{a}dl\`{a}g control path. We therefore, introduce the
measurable projection that forgets the latent refresh labels and retains only
the realised control trajectory. Define
\(
  \mathfrak R_{t,T}:\widehat{\mathcal U}_{t,T}\to\mathcal U_{t,T}
\)
as follows. If \(N=0\), then \(\mathfrak R_{t,T}(u,0)\) is the constant path
equal to \(u\) on \([t,T]\). If
\[
  \mathfrak m=\sum_{k=1}^N\delta_{(\tau_k,v_k)},
  \qquad
  t<\tau_1<\cdots<\tau_N\le T,
  \qquad N\ge 1,
\]
then \(\mathfrak R_{t,T}(u,\mathfrak m)\) is the c\`{a}dl\`{a}g path
\[
  \mathfrak R_{t,T}(u,\mathfrak m)(s)
  =
  \begin{cases}
    u, & s\in[t,\tau_1),\\
    v_k, & s\in[\tau_k,\tau_{k+1}),\quad k=1,\ldots,N-1,\\
    v_N, & s\in[\tau_N,T].
  \end{cases}
\]
The map \(\mathfrak R_{t,T}\) is Borel measurable.

\begin{corollary}[Observed control-path inequality]
\label{cor:observed-control-path-inequality}
Let
\[
  \mathbf Q^{\alpha;\mathbf x,u}
  :=
  (\mathfrak R_{t,T})_\#
  \widehat{\mathbf Q}^{\alpha;\mathbf x,u},
  \qquad
  \mathbf Q^{I;\mathbf x,u}
  :=
  (\mathfrak R_{t,T})_\#
  \widehat{\mathbf Q}^{I;\mathbf x,u}
\]
be the corresponding conditional laws on the observed control-path space
\(\mathcal U_{t,T}\). Then, for every
\((\mathbf x,u)\in\mathcal X_{t,T}\times\mathbb R^q\),
\begin{equation}
  D\!\Bigl(
    \mathbf Q^{\alpha;\mathbf x,u}
    \,\Big\|\,
    \mathbf Q^{I;\mathbf x,u}
  \Bigr)
  \le
  \lambda
  \int_t^T
    D\!\Bigl(
      c_s^\alpha(\cdot\mid \mathbf x(s),\mu_s^{t,\zeta,\alpha})
      \,\Big\|\,
      c_s^I(\cdot\mid \mathbf x(s),\mu_s^{t,\zeta,\alpha})
    \Bigr)\,
    \mathrm ds.
  \label{eq:observed-control-path-inequality}
\end{equation}
\end{corollary}

\begin{proof}
This follows from Proposition~\ref{prop:control-KL-poisson-augmented} and the
data-processing inequality for relative entropy applied to the measurable map
\(\mathfrak R_{t,T}\) \cite{DemboZeitouni1998LDP}. The inequality need not be an equality, since \(\mathfrak R_{t,T}\) forgets the latent refresh labels; in particular, a
refresh may redraw the current control value and leave the observed control path
unchanged.
\end{proof}

\section{Reference drift and local state cost}
\label{sec:state-cost}

We next introduce the state-side component of the running cost. The purpose of
this term is to measure, locally in time, how far the realised state dynamics
are from a prescribed reference dynamics. To define this reference dynamics, fix
a reference drift
\(
  f^I:[0,T]\times\mathbb R^n\times\mathbb R^q\to\mathbb R^n .
\)
The reference drift is evaluated at the same state and control variables as the
actual drift \(f\).

\begin{assumption}\label{assump:reference-drift-poisson}
The following hold.
\begin{enumerate}
\item
The map \(f^I\) is jointly Borel measurable, continuous in \(t\), and globally
Lipschitz in \((x,u)\), uniformly in \(t\).

\item
There exists a constant \(C_I>0\) such that
\(
  \|f^I(t,x,u)\|
  \le
  C_I\bigl(1+\|x\|+\|u\|\bigr)
\)
for all \((t,x,u)\in[0,T]\times\mathbb R^n\times\mathbb R^q\).
\end{enumerate}
\end{assumption}

Recall from Section~\ref{sec:poisson-refresh-formulation} that
\(
  \Sigma(x,u)
  :=
  \sigma(x,u)\sigma(x,u)^\top\), for every
  \((x,u)\in\mathbb R^n\times\mathbb R^q
\).
By Assumption~\ref{assump:primitive-coeff}(2), the matrix \(\Sigma(x,u)\) is
uniformly positive definite, hence invertible, for every
\((x,u)\in\mathbb R^n\times\mathbb R^q\). We therefore define the local
state-side cost density by
\begin{equation}
  L(t,x,u)
  :=
  \frac12\,
  \bigl(
    f(x,u)-f^I(t,x,u)
  \bigr)^\top
  \Sigma(x,u)^{-1}
  \bigl(
    f(x,u)-f^I(t,x,u)
  \bigr),
  \label{eq:local-state-integrand-poisson}
\end{equation}
for
\(
  (t,x,u)\in[0,T]\times\mathbb R^n\times\mathbb R^q 
\).
Thus \(L(t,x,u)\) measures the local drift mismatch between the actual drift
\(f(x,u)\) and the reference drift \(f^I(t,x,u)\), weighted by the inverse
diffusion covariance. It is evaluated at the realised control value \(u\),
because in the Poisson-refresh formulation the state equation is driven by the
actual current control value, not by an average over the refresh kernel.

The following lemma records the basic measurability, nonnegativity, and growth
properties of the pointwise integrand \(L\).

\begin{lemma}\label{lem:state-cost-growth-poisson}
Under Assumptions~\ref{assump:primitive-coeff} and
\ref{assump:reference-drift-poisson}, the function \(L\) is Borel measurable
and nonnegative. Moreover, there exists a constant \(C>0\), depending only on
the constants in Assumptions~\ref{assump:primitive-coeff} and
\ref{assump:reference-drift-poisson}, such that
\begin{equation}
  0\le L(t,x,u)\le C\bigl(1+\|x\|^2+\|u\|^2\bigr)
  \label{eq:L-growth-poisson}
\end{equation}
for all \((t,x,u)\in[0,T]\times\mathbb R^n\times\mathbb R^q\).
\end{lemma}

\begin{proof}
Since \(\sigma\) is continuous and
\(\Sigma=\sigma\sigma^\top\) is uniformly positive definite by
\eqref{eq:uniform-ellipticity-poisson}, the map
\(
  (x,u)\longmapsto \Sigma(x,u)^{-1}
\)
is continuous and uniformly bounded on \(\mathbb R^n\times\mathbb R^q\). Hence
\(L\) is Borel measurable. Its nonnegativity follows from the positive
definiteness of \(\Sigma(x,u)^{-1}\).

Finally, since \(\Sigma^{-1}\) is uniformly bounded and both \(f\) and \(f^I\)
have linear growth, there exists \(C>0\) such that
\[
  \|f(x,u)-f^I(t,x,u)\|^2
  \le
  C\bigl(1+\|x\|^2+\|u\|^2\bigr).
\]
Substituting this estimate into the definition of \(L\) gives
\eqref{eq:L-growth-poisson}.
\end{proof}

We now average the pointwise integrand \(L\) against a joint law of the current
state and control. For
\(
  \eta\in\mathcal P_2(\mathbb R^n\times\mathbb R^q),
\)
define the local state functional
\begin{equation}
  \mathcal J_{\mathrm{state}}(t,\eta)
  :=
  \int_{\mathbb R^n\times\mathbb R^q}
    L(t,x,u)\,\eta(\mathrm dx,\mathrm du),
  \qquad
  (t,\eta)\in
  [0,T]\times\mathcal P_2(\mathbb R^n\times\mathbb R^q).
  \label{eq:local-state-functional-poisson}
\end{equation}
The growth bound in Lemma~\ref{lem:state-cost-growth-poisson} ensures that this
quantity is finite for every
\((t,\eta)\in[0,T]\times\mathcal P_2(\mathbb R^n\times\mathbb R^q)\). Indeed,
integrating \eqref{eq:L-growth-poisson} against \(\eta\) gives the estimate
recorded in the next lemma. We also record measurability, since this functional
will later be integrated along closed-loop law flows.

\begin{lemma}\label{lem:local-state-functional-measurable}
The map
\(
  (t,\eta)\longmapsto \mathcal J_{\mathrm{state}}(t,\eta)
\)
is Borel measurable on
\(
  [0,T]\times\mathcal P_2(\mathbb R^n\times\mathbb R^q),
\)
and there exists a constant \(C>0\) such that
\begin{equation}
  0\le \mathcal J_{\mathrm{state}}(t,\eta)
  \le
  C\bigl(1+m_2(\eta)\bigr)
  \label{eq:local-state-functional-growth}
\end{equation}
for all
\(
  (t,\eta)\in[0,T]\times\mathcal P_2(\mathbb R^n\times\mathbb R^q),
\)
where \(m_2(\eta)\) denotes the second moment of \(\eta\) on
\(\mathbb R^n\times\mathbb R^q\). In particular,
\(
  \mathcal J_{\mathrm{state}}(t,\eta)<\infty
\)
for every such \((t,\eta)\).
\end{lemma}

\begin{proof}
By Lemma~\ref{lem:state-cost-growth-poisson},
\[
  0\le L(t,x,u)\le C\bigl(1+\|x\|^2+\|u\|^2\bigr).
\]
Integrating this inequality with respect to \(\eta\) gives
\eqref{eq:local-state-functional-growth}.

It remains to prove measurability. For \(N\in\mathbb N\), define
\(
  L_N(t,x,u):=L(t,x,u)\wedge N.
\)
Each \(L_N\) is bounded and Borel measurable. Hence
\[
  (t,\eta)\longmapsto
  \int_{\mathbb R^n\times\mathbb R^q}
    L_N(t,x,u)\,\eta(\mathrm dx,\mathrm du)
\]
is Borel measurable. Since
\(
  L_N(t,x,u)\uparrow L(t,x,u)\) as \(N\to\infty \),
the monotone convergence theorem gives
\[
  \mathcal J_{\mathrm{state}}(t,\eta)
  =
  \sup_{N\in\mathbb N}
  \int_{\mathbb R^n\times\mathbb R^q}
    L_N(t,x,u)\,\eta(\mathrm dx,\mathrm du).
\]
Therefore \((t,\eta)\mapsto\mathcal J_{\mathrm{state}}(t,\eta)\) is Borel
measurable as the pointwise supremum of Borel measurable functions.
\end{proof}

We now evaluate the local state functional along the closed-loop law generated
by an admissible sampling rule. Since the integrand \(L(t,x,u)\) depends on the
realised control value \(u\), it is evaluated against the joint law of
\((X_s,U_s)\). For \(0\le t\le T\), \(\zeta\in\mathcal P_2(\mathbb R^n\times\mathbb R^q)\),
and \(\alpha\in\mathfrak A[t,T]\), let
\(
  \eta_s^{t,\zeta,\alpha}
  :=
  \mathcal L_{\mathbb P^{t,\zeta,\alpha}}(X_s,U_s)\), for every \(s\in[t,T]\)
be the joint law flow from
Theorem~\ref{thm:wellposedness-poisson-refresh}. The state contribution
of \(\alpha\) on \([t,T]\) is defined by
\begin{align}
  \mathcal H_{\mathrm{state}}(t,\zeta;\alpha)
  &:=
  \int_t^T
    \mathcal J_{\mathrm{state}}\!\bigl(
      s,\eta_s^{t,\zeta,\alpha}
    \bigr)\,\mathrm ds
  \label{eq:state-cost-general-time-poisson}
  \\
  &=
  \mathbb E^{t,\zeta,\alpha}\!\left[
    \int_t^T
      L(s,X_s,U_s)\,\mathrm ds
  \right],
  \label{eq:state-path-cost-poisson}
\end{align}
with values in \([0,\infty]\). The equality of
\eqref{eq:state-cost-general-time-poisson} and
\eqref{eq:state-path-cost-poisson} follows by expanding the definition of
\(\mathcal J_{\mathrm{state}}\) and using
\(\eta_s^{t,\zeta,\alpha}
=
\mathcal L_{\mathbb P^{t,\zeta,\alpha}}(X_s,U_s)\).

The next proposition records that this state contribution is always finite
under the standing assumptions.

\begin{proposition}\label{prop:state-cost-finite-poisson}
For every \(0\le t\le T\), every
\(\zeta\in\mathcal P_2(\mathbb R^n\times\mathbb R^q)\), and every
\(\alpha\in\mathfrak A[t,T]\), one has
\(
  \mathcal H_{\mathrm{state}}(t,\zeta;\alpha)<\infty
\).
More precisely, there exists a constant \(C>0\), depending only on the
constants in Assumptions~\ref{assump:primitive-coeff} and
\ref{assump:reference-drift-poisson}, on \(L_\alpha,M_\alpha\), and on \(T\),
such that
\begin{equation}
  \mathcal H_{\mathrm{state}}(t,\zeta;\alpha)
  \le
  C\bigl(1+m_2(\zeta)\bigr).
  \label{eq:state-cost-bound-poisson}
\end{equation}
\end{proposition}

\begin{proof}
Fix \(0\le t\le T\), \(\zeta\in\mathcal P_2(\mathbb R^n\times\mathbb R^q)\),
and \(\alpha\in\mathfrak A[t,T]\). By \eqref{eq:L-growth-poisson},
\[
  L(s,X_s,U_s)
  \le
  C\bigl(
    1+\|X_s\|^2+\|U_s\|^2
  \bigr).
\]
Taking expectation, integrating over \(s\in[t,T]\), and using
\eqref{eq:moment-bound-joint-poisson} from
Theorem~\ref{thm:wellposedness-poisson-refresh}, we obtain
\begin{align*}
  \mathcal H_{\mathrm{state}}(t,\zeta;\alpha)
  &=
  \mathbb E^{t,\zeta,\alpha}\!\left[
    \int_t^T
      L(s,X_s,U_s)\,\mathrm ds
  \right]
  \\
  &\le
  C\int_t^T
    \Bigl(
      1+
      \mathbb E^{t,\zeta,\alpha}[\|X_s\|^2+\|U_s\|^2]
    \Bigr)\,\mathrm ds
  \\
  &\le
  C\bigl(1+m_2(\zeta)\bigr).
\end{align*}
This proves \eqref{eq:state-cost-bound-poisson}.
\end{proof}

The choice of the integrand \(L\) is motivated by the following path-space
identity for diffusions driven by a fixed deterministic control path \cite{KaratzasShreve1991BMSC,LiptserShiryaev2001SRP}.

For \(0\le t\le T\), let \( \mathcal U_2[t,T]:=
  L^2([t,T];\mathbb R^q) \)
denote the space of deterministic square-integrable control paths. Fix
\((t,x,\mathbf u)\in[0,T]\times\mathbb R^n\times\mathcal U_2[t,T]\), and
consider the two diffusions driven by the same control path \(\mathbf u\):
\begin{align}
  \mathrm dX_s
  &=
  f\bigl(X_s,\mathbf u(s)\bigr)\,\mathrm ds
  +
  \sigma\bigl(X_s,\mathbf u(s)\bigr)\,\mathrm dW_s,
  \qquad
  X_t=x,
  \label{eq:path-controlled-SDE-poisson}
  \\
  \mathrm dX_s^I
  &=
  f^I\bigl(s,X_s^I,\mathbf u(s)\bigr)\,\mathrm ds
  +
  \sigma\bigl(X_s^I,\mathbf u(s)\bigr)\,\mathrm d\widetilde W_s,
  \qquad
  X_t^I=x.
  \label{eq:path-reference-SDE-poisson}
\end{align}
Under Assumptions~\ref{assump:primitive-coeff} and
\ref{assump:reference-drift-poisson}, each of these equations has a unique
strong solution. Indeed, for fixed \(\mathbf u\in\mathcal U_2[t,T]\), the
coefficients in \eqref{eq:path-controlled-SDE-poisson} and
\eqref{eq:path-reference-SDE-poisson} are Borel in time, globally Lipschitz in
\(x\) uniformly in time, and satisfy the standard integrability bounds because
\(\mathbf u\in L^2([t,T];\mathbb R^q)\). Let
\(
  \mathbf P^{t,x,\mathbf u}
  :=
  \mathcal L(X) \), and \(
  \mathbf P^{I,t,x,\mathbf u}
  :=
  \mathcal L(X^I) \)
denote the corresponding path laws on \(
  \mathcal X_{t,T}=C([t,T];\mathbb R^n) \)
The following proposition identifies the relative entropy between these two
path laws.

\begin{proposition}[Frozen-control entropy identity]
\label{prop:pathwise-state-KL-poisson}
Fix \((t,x,\mathbf u)\in[0,T]\times\mathbb R^n\times\mathcal U_2[t,T]\).
Assume that
\begin{equation}
  \mathbb E_{\mathbf P^{I,t,x,\mathbf u}}\!\left[
    \exp\!\left(
      \int_t^T
        L\bigl(s,X_s,\mathbf u(s)\bigr)\,\mathrm ds
    \right)
  \right]
  <\infty.
  \label{eq:Novikov-conditional-poisson}
\end{equation}
Then
\begin{equation}
  D\!\Bigl(
    \mathbf P^{t,x,\mathbf u}
    \,\Big\|\,
    \mathbf P^{I,t,x,\mathbf u}
  \Bigr)
  =
  \mathbb E_{\mathbf P^{t,x,\mathbf u}}\!\left[
    \int_t^T
      L\bigl(s,X_s,\mathbf u(s)\bigr)\,\mathrm ds
  \right].
  \label{eq:pathwise-state-KL-conditional-poisson}
\end{equation}
\end{proposition}

\begin{proof}
Work on the canonical path space \(\mathcal X_{t,T}\), and let \(X\) denote the
canonical process. Write
\[
  \delta_s
  :=
  f\bigl(X_s,\mathbf u(s)\bigr)
  -
  f^I\bigl(s,X_s,\mathbf u(s)\bigr),
  \qquad
  \Sigma_s
  :=
  \Sigma\bigl(X_s,\mathbf u(s)\bigr).
\]
By Assumption~\ref{assump:primitive-coeff}, \(\Sigma_s\) is uniformly elliptic,
so \(\Sigma_s^{-1}\) is well defined and uniformly bounded. 

Under \(\mathbf P^{I,t,x,\mathbf u}\), the canonical process solves the
reference diffusion. Hence
\[
  M_s
  :=
  X_s-x-\int_t^s
    f^I\bigl(r,X_r,\mathbf u(r)\bigr)\,\mathrm dr
\]
is a continuous local martingale with quadratic variation
\[
  \langle M\rangle_s
  =
  \int_t^s \Sigma_r\,\mathrm dr .
\]
By uniform ellipticity, \(\Sigma_r^{-1}\) is well defined and uniformly
bounded. Therefore the stochastic integral
\[
  \int_t^T
    \delta_s^\top\Sigma_s^{-1}\,\mathrm dM_s
\]
is well defined. Define
\[
  Z_T
  :=
  \exp\!\left(
    \int_t^T
      \delta_s^\top\Sigma_s^{-1}\,\mathrm dM_s
    -
    \frac12\int_t^T
      \delta_s^\top\Sigma_s^{-1}\delta_s\,\mathrm ds
  \right).
\]
Since
\[
  L\bigl(s,X_s,\mathbf u(s)\bigr)
  =
  \frac12
  \delta_s^\top\Sigma_s^{-1}\delta_s,
\]
assumption \eqref{eq:Novikov-conditional-poisson} is exactly Novikov's
condition for \(Z_T\). Thus \(Z_T\) is a martingale with expectation one, and we
may define a probability measure \(\mathbf Q\) on \(\mathcal X_{t,T}\) by
\[
  \frac{\mathrm d\mathbf Q}
       {\mathrm d\mathbf P^{I,t,x,\mathbf u}}
  :=
  Z_T.
\]
By Girsanov's theorem for continuous martingales \cite{KaratzasShreve1991BMSC,LiptserShiryaev2001SRP}, under \(\mathbf Q\) the
process
\[
  M_s-\int_t^s \delta_r\,\mathrm dr
\]
is a continuous local martingale with quadratic variation
\(\langle M\rangle_s\). Consequently, under \(\mathbf Q\), the
canonical process \(X\) solves
\[
  \mathrm dX_s
  =
  f\bigl(X_s,\mathbf u(s)\bigr)\,\mathrm ds
  +
  \sigma\bigl(X_s,\mathbf u(s)\bigr)\,\mathrm dW_s
\]
for some Brownian motion \(W\). By uniqueness in law for
\eqref{eq:path-controlled-SDE-poisson}, we have
\(
  \mathbf Q=\mathbf P^{t,x,\mathbf u}
\).

Finally, we compute the relative entropy. Since
\(\mathrm d\mathbf P^{t,x,\mathbf u}/\mathrm d\mathbf P^{I,t,x,\mathbf u}=Z_T\),
\[
  D\!\Bigl(
    \mathbf P^{t,x,\mathbf u}
    \,\Big\|\,
    \mathbf P^{I,t,x,\mathbf u}
  \Bigr)
  =
  \mathbb E_{\mathbf P^{t,x,\mathbf u}}[\log Z_T].
\]
Under \(\mathbf P^{t,x,\mathbf u}\),
\[
  M_s-\int_t^s\delta_r\,\mathrm dr
\]
is a continuous local martingale. Localizing this martingale and then passing to
the limit gives
\[
  \mathbb E_{\mathbf P^{t,x,\mathbf u}}[\log Z_T]
  =
  \frac12
  \mathbb E_{\mathbf P^{t,x,\mathbf u}}\!\left[
    \int_t^T
      \delta_s^\top\Sigma_s^{-1}\delta_s\,\mathrm ds
  \right].
\]
Since the integrand is \(2L(s,X_s,\mathbf u(s))\), this proves
\eqref{eq:pathwise-state-KL-conditional-poisson}.
\end{proof}

We now combine the control-side and state-side local costs. Let
\[
  \operatorname{pr}_x:\mathbb R^n\times\mathbb R^q\to\mathbb R^n,
  \qquad
  \operatorname{pr}_x(x,u):=x,
\]
denote the first-coordinate projection. For
\(t\in[0,T]\), \(\eta\in\mathcal P_2(\mathbb R^n\times\mathbb R^q)\), and
\(\gamma\in\mathcal K\), define
\begin{equation}
  \mathcal J(t,\eta;\gamma)
  :=
  \mathcal J_{\mathrm{ctrl}}\!\bigl(
    t,(\operatorname{pr}_x)_\#\eta;\gamma
  \bigr)
  +
  \mathcal J_{\mathrm{state}}(t,\eta)
  \in[0,\infty].
  \label{eq:J-inst-general-poisson}
\end{equation}
Here \((\operatorname{pr}_x)_\#\eta\) is the state marginal of the joint law
\(\eta\).

For \(0\le t\le T\), \(\zeta\in\mathcal P_2(\mathbb R^n\times\mathbb R^q)\),
and \(\alpha\in\mathfrak A[t,T]\), define the finite-horizon cost by
\begin{equation}
  \mathcal J_{t,T}(\zeta,\alpha)
  :=
  \int_t^T
    \mathcal J\!\bigl(
      s,\eta_s^{t,\zeta,\alpha};\gamma_s^{t,\zeta,\alpha}
    \bigr)\,\mathrm ds
  =
  \mathcal H_{\mathrm{ctrl}}(t,\zeta;\alpha)
  +
  \mathcal H_{\mathrm{state}}(t,\zeta;\alpha),
  \label{eq:value-cost-functional-local-poisson}
\end{equation}
with values in \([0,\infty]\), where
\(\gamma_s^{t,\zeta,\alpha}\) is the instantaneous closed-loop refresh kernel
defined in \eqref{eq:gamma-from-alpha-poisson}. Since
\(\mathcal H_{\mathrm{state}}(t,\zeta;\alpha)<\infty\) by
Proposition~\ref{prop:state-cost-finite-poisson} and
\(\mathcal H_{\mathrm{ctrl}}(t,\zeta;\alpha)\in[0,\infty]\) by definition,
\(\mathcal J_{t,T}(\zeta,\alpha)\) is well defined as an element of
\([0,\infty]\).

The initial optimization problem is
\begin{equation}
  \inf_{\alpha\in\mathfrak A}\mathcal J_{0,T}(\zeta_0,\alpha),
  \label{eq:initial-optimization-problem-poisson}
\end{equation}
where \(\zeta_0\in\mathcal P_2(\mathbb R^n\times\mathbb R^q)\) is the prescribed
initial joint law. If one starts instead from a prescribed initial state law
\(\mu_0\in\mathcal P_2(\mathbb R^n)\), then an initial control law must also be
specified. Equivalently, one chooses a kernel
\(
  \rho_0(\cdot\mid x,\mu_0)\in\mathcal P(\mathbb R^q),
\)
and sets
\[
  \zeta_0(\mathrm dx,\mathrm du)
  :=
  \mu_0(\mathrm dx)\,\rho_0(\mathrm du\mid x,\mu_0).
\]

If one also wishes to benchmark the prescribed initial joint law \(\zeta_0\)
against a reference initial joint law
\(\zeta_0^I\in\mathcal P(\mathbb R^n\times\mathbb R^q)\) with
\[
  D(\zeta_0\|\zeta_0^I)<\infty,
\]
then one may add the constant term \(D(\zeta_0\|\zeta_0^I)\) to
\eqref{eq:initial-optimization-problem-poisson}. Since \(\zeta_0\) is fixed throughout, this additive constant does not affect the minimization or the dynamic programming principle. We therefore omit it in the continuation problem considered in the next section.

\section{Continuation problem and dynamic programming principle}
\label{sec:dpp-poisson}
We now formulate the continuation problem obtained by restarting
\eqref{eq:initial-optimization-problem-poisson} at an arbitrary time \(t\) with
an arbitrary joint law \(\zeta\) of \((X_t,U_t)\). The joint law is the natural
state variable for continuation, since the current control value \(U_t\)
persists until the next refresh time.

Recall from Sections~\ref{sec:control-entropy} and \ref{sec:state-cost} that
for \(0\le t\le T\), \(\zeta\in\mathcal P_2(\mathbb R^n\times\mathbb R^q)\),
and \(\alpha\in\mathfrak A[t,T]\),
\(
  \mathbb P^{t,\zeta,\alpha}\in\mathcal P(\mathcal Y_{t,T})
\)
denotes the closed-loop law, with joint law flow
\[
  \eta_s^{t,\zeta,\alpha}
  :=
  \mathcal L_{\mathbb P^{t,\zeta,\alpha}}(X_s,U_s),
  \qquad
  s\in[t,T],
\]
and state-law flow
\[
  \mu_s^{t,\zeta,\alpha}
  :=
  \mathcal L_{\mathbb P^{t,\zeta,\alpha}}(X_s)
  =
  (\operatorname{pr}_x)_\#\eta_s^{t,\zeta,\alpha}.
\]
The corresponding instantaneous refresh kernel is
\[
  \gamma_s^{t,\zeta,\alpha}(\cdot\mid x)
  :=
  c_s^\alpha(\cdot\mid x,\mu_s^{t,\zeta,\alpha}),
  \qquad
  s\in[t,T],\ x\in\mathbb R^n.
\]

To split costs at intermediate times, we use the same running-cost functional on
subintervals. Thus, for
\(0\le t\le \tau\le T\), \(\zeta\in\mathcal P_2(\mathbb R^n\times\mathbb R^q)\),
and \(\alpha\in\mathfrak A[t,\tau]\), define
\begin{equation}
  \mathcal J_{t,\tau}(\zeta,\alpha)
  :=
  \int_t^\tau
    \mathcal J\!\bigl(
      s,\eta_s^{t,\zeta,\alpha};\gamma_s^{t,\zeta,\alpha}
    \bigr)\,\mathrm ds
  \in [0,\infty],
  \label{eq:value-cost-functional-short-poisson}
\end{equation}
where \(\mathcal J\) is the local running-cost functional from
\eqref{eq:J-inst-general-poisson}. This is the restriction of the
finite-horizon cost in \eqref{eq:value-cost-functional-local-poisson} to the
shorter interval \([t,\tau]\). Since the integrand is nonnegative,
\(\mathcal J_{t,\tau}(\zeta,\alpha)\) is well defined as an element of
\([0,\infty]\).

The continuation value is the optimal cost from the current time \(t\) all the
way to the terminal time \(T\):
\begin{equation}
  V(t,\zeta)
  :=
  \inf_{\alpha\in\mathfrak A[t,T]}
  \mathcal J_{t,T}(\zeta,\alpha),
  \qquad
  (t,\zeta)\in[0,T]\times\mathcal P_2(\mathbb R^n\times\mathbb R^q),
  \label{eq:value-general-poisson}
\end{equation}
with terminal condition \(V(T,\zeta)=0 \)

As observed in Section~\ref{sec:state-cost}, if an initial relative-entropy term \(D(\zeta_0\|\zeta_0^I)\) is also included,
then \(D(\zeta_0\|\zeta_0^I)+\mathcal J_{0,T}(\zeta_0,\alpha) \)
differs from \(\mathcal J_{0,T}(\zeta_0,\alpha)\) by a constant independent of
\(\alpha\). Therefore this term does not affect the minimization or the
continuation value function in \eqref{eq:value-general-poisson}.

The next proposition records that \(V\) is finite and has at most quadratic
growth in the initial joint law.

\begin{proposition}\label{prop:value-finite-poisson}
There exists a constant \(C>0\), depending only on the constants in
Assumptions~\ref{assump:primitive-coeff} and
\ref{assump:reference-drift-poisson}, on the admissibility constants of the
reference sampling rule \(\alpha^I\), and on \(T\), such that
\begin{equation}
  0\le V(t,\zeta)\le C\bigl(1+m_2(\zeta)\bigr),
  \qquad
  (t,\zeta)\in[0,T]\times\mathcal P_2(\mathbb R^n\times\mathbb R^q).
  \label{eq:value-growth-bound-poisson}
\end{equation}
\end{proposition}

\begin{proof}
Nonnegativity is immediate from the definition of \(\mathcal J\). For the upper
bound, use the admissible reference sampling rule \(\alpha^I|_{[t,T]}\in\mathfrak A[t,T]\). 
By \eqref{eq:local-control-zero-reference-poisson},
\(
  \mathcal J_{\mathrm{ctrl}}\!\bigl(
    s,\mu_s^{t,\zeta,\alpha^I};\gamma_s^{t,\zeta,\alpha^I}
  \bigr)
  =0
\)
for every \(s\in[t,T]\). Therefore
\[
  \mathcal J_{t,T}\!\bigl(
    \zeta,\alpha^I|_{[t,T]}
  \bigr)
  =
  \mathcal H_{\mathrm{state}}\!\bigl(
    t,\zeta;\alpha^I|_{[t,T]}
  \bigr).
\]
Proposition~\ref{prop:state-cost-finite-poisson} gives
\[
  \mathcal J_{t,T}\!\bigl(
    \zeta,\alpha^I|_{[t,T]}
  \bigr)
  \le
  C\bigl(1+m_2(\zeta)\bigr).
\]
Taking the infimum over \(\alpha\in\mathfrak A[t,T]\) proves
\eqref{eq:value-growth-bound-poisson}.
\end{proof}

The dynamic programming principle relies on decomposing the cost over
subintervals. We use the restriction and concatenation operations introduced in
Proposition~\ref{prop:stability-sampling-rules}: for
\(0\le t\le\tau\le T\), an admissible rule on \([t,T]\) can be restricted to
\([t,\tau]\) and \([\tau,T]\), while
\(\alpha^1\in\mathfrak A[t,\tau]\) and
\(\alpha^2\in\mathfrak A[\tau,T]\) can be concatenated to form
\(
  \alpha^1\otimes_\tau\alpha^2\in\mathfrak A[t,T]
\).

The following lemma gives the corresponding decomposition of the cost
functional.

\begin{lemma}[Cost splitting under restriction and concatenation]
\label{lem:cost-splitting-poisson}
Let \(0\le t\le\tau\le T\).

\begin{enumerate}
\item
If \(\alpha\in\mathfrak A[t,T]\), then
\begin{equation}
  \mathcal J_{t,T}(\zeta,\alpha)
  =
  \mathcal J_{t,\tau}\!\bigl(
    \zeta,\alpha|_{[t,\tau]}
  \bigr)
  +
  \mathcal J_{\tau,T}\!\bigl(
    \eta_\tau^{t,\zeta,\alpha},
    \alpha|_{[\tau,T]}
  \bigr).
  \label{eq:cost-splitting-restriction-poisson}
\end{equation}

\item
If \(\alpha^1\in\mathfrak A[t,\tau]\) and
\(\alpha^2\in\mathfrak A[\tau,T]\), and if
\(
  \vartheta:=\eta_\tau^{t,\zeta,\alpha^1},
\)
then
\begin{equation}
  \mathcal J_{t,T}\!\bigl(
    \zeta,\alpha^1\otimes_\tau\alpha^2
  \bigr)
  =
  \mathcal J_{t,\tau}(\zeta,\alpha^1)
  +
  \mathcal J_{\tau,T}(\vartheta,\alpha^2).
  \label{eq:cost-splitting-concatenation-poisson}
\end{equation}
\end{enumerate}
\end{lemma}

\begin{proof}
For part~(1), Proposition~\ref{prop:stability-joint-laws}(1) gives
\(
  \eta_s^{t,\zeta,\alpha}
  =
  \eta_s^{t,\zeta,\alpha|_{[t,\tau]}}\) for every \(s\in[t,\tau]\).
Taking first marginals gives the corresponding identity for the state-law
flows. Since the sampling rule on \([t,\tau]\) is
\(\alpha|_{[t,\tau]}\), the associated refresh kernels also agree on
\([t,\tau]\).

For the second subinterval, apply
Proposition~\ref{prop:stability-joint-laws}(2) with \(\alpha^1:=\alpha|_{[t,\tau]}\), and \( \alpha^2:=\alpha|_{[\tau,T]}\).
Since \( \alpha^1\otimes_\tau\alpha^2=\alpha\),
we obtain
\[
  \eta_s^{t,\zeta,\alpha}
  =
  \eta_s^{\tau,\eta_\tau^{t,\zeta,\alpha},\,\alpha|_{[\tau,T]}},
  \qquad
  s\in[\tau,T].
\]
Again, taking first marginals gives the corresponding identity for the
state-law flows, and the restricted sampling rules imply equality of the
associated refresh kernels on \([\tau,T]\). Splitting the integral in
\eqref{eq:value-cost-functional-short-poisson} at time \(\tau\) gives
\eqref{eq:cost-splitting-restriction-poisson}.

For part~(2), Proposition~\ref{prop:stability-joint-laws}(2) gives
\[
  \eta_s^{t,\zeta,\alpha^1\otimes_\tau\alpha^2}
  =
  \begin{cases}
    \eta_s^{t,\zeta,\alpha^1}, & s\in[t,\tau],\\
    \eta_s^{\tau,\vartheta,\alpha^2}, & s\in[\tau,T].
  \end{cases}
\]
Taking first marginals gives the same piecewise identity for the state-law
flows. Since the concatenated sampling rule agrees with \(\alpha^1\) on
\([t,\tau]\) and with \(\alpha^2\) on \((\tau,T]\), the associated refresh
kernels agree on the corresponding subintervals. Splitting the defining integral for
\(\mathcal J_{t,T}\) in
\eqref{eq:value-cost-functional-short-poisson} at \(\tau\) yields
\eqref{eq:cost-splitting-concatenation-poisson}.
\end{proof}

We can now prove the dynamic programming principle \cite{DjetePossamaiTan2022MKVDPP}.

\begin{theorem}[Dynamic programming principle]
\label{thm:DPP-poisson}
Fix \(0\le t\le\tau\le T\) and
\(\zeta\in\mathcal P_2(\mathbb R^n\times\mathbb R^q)\). Then
\begin{equation}
  V(t,\zeta)
  =
  \inf_{\alpha\in\mathfrak A[t,\tau]}
  \left\{
    \mathcal J_{t,\tau}(\zeta,\alpha)
    +
    V\!\bigl(
      \tau,\eta_\tau^{t,\zeta,\alpha}
    \bigr)
  \right\}.
  \label{eq:DPP-general-poisson}
\end{equation}
Equivalently,
\begin{equation}
  V(t,\zeta)
  =
  \inf_{\alpha\in\mathfrak A[t,\tau]}
  \left\{
    \int_t^\tau
      \mathcal J\!\bigl(
        s,\eta_s^{t,\zeta,\alpha};\gamma_s^{t,\zeta,\alpha}
      \bigr)\,\mathrm ds
    +
    V\!\bigl(
      \tau,\eta_\tau^{t,\zeta,\alpha}
    \bigr)
  \right\}.
  \label{eq:DPP-general-J-poisson}
\end{equation}
\end{theorem}

\begin{proof}
Set
\[
  \Xi(t,\tau,\zeta)
  :=
  \inf_{\alpha\in\mathfrak A[t,\tau]}
  \left\{
    \mathcal J_{t,\tau}(\zeta,\alpha)
    +
    V\!\bigl(
      \tau,\eta_\tau^{t,\zeta,\alpha}
    \bigr)
  \right\}.
\]

We first show that
\(
  V(t,\zeta)\ge \Xi(t,\tau,\zeta)
\).
Let \(\alpha\in\mathfrak A[t,T]\). By
\eqref{eq:cost-splitting-restriction-poisson},
\[
  \mathcal J_{t,T}(\zeta,\alpha)
  =
  \mathcal J_{t,\tau}\!\bigl(
    \zeta,\alpha|_{[t,\tau]}
  \bigr)
  +
  \mathcal J_{\tau,T}\!\bigl(
    \eta_\tau^{t,\zeta,\alpha},
    \alpha|_{[\tau,T]}
  \bigr).
\]
By definition of the value function,
\(
  V\!\bigl(
    \tau,\eta_\tau^{t,\zeta,\alpha}
  \bigr)
  \le
  \mathcal J_{\tau,T}\!\bigl(
    \eta_\tau^{t,\zeta,\alpha},
    \alpha|_{[\tau,T]}
  \bigr) \).
Therefore
\[
  \mathcal J_{t,T}(\zeta,\alpha)
  \ge
  \mathcal J_{t,\tau}\!\bigl(
    \zeta,\alpha|_{[t,\tau]}
  \bigr)
  +
  V\!\bigl(
    \tau,\eta_\tau^{t,\zeta,\alpha}
  \bigr)
  \ge
  \Xi(t,\tau,\zeta).
\]
Taking the infimum over \(\alpha\in\mathfrak A[t,T]\) yields \( V(t,\zeta)\ge \Xi(t,\tau,\zeta)\).

We next show that
\(
  V(t,\zeta)\le \Xi(t,\tau,\zeta)
\).
Fix \(\varepsilon>0\). By the definition of infimum, there exists
\(\alpha_\varepsilon^1\in\mathfrak A[t,\tau]\) such that
\[
  \mathcal J_{t,\tau}(\zeta,\alpha_\varepsilon^1)
  +
  V\!\bigl(
    \tau,\eta_\tau^{t,\zeta,\alpha_\varepsilon^1}
  \bigr)
  \le
  \Xi(t,\tau,\zeta)+\varepsilon.
\]
Set
\(
  \vartheta_\varepsilon
  :=
  \eta_\tau^{t,\zeta,\alpha_\varepsilon^1}.
\)
By Proposition~\ref{prop:value-finite-poisson},
\(V(\tau,\vartheta_\varepsilon)<\infty\). Again by the definition of infimum,
there exists \(\alpha_\varepsilon^2\in\mathfrak A[\tau,T]\) such that
\[
  \mathcal J_{\tau,T}(\vartheta_\varepsilon,\alpha_\varepsilon^2)
  \le
  V(\tau,\vartheta_\varepsilon)+\varepsilon.
\]
Let
\(
  \alpha_\varepsilon
  :=
  \alpha_\varepsilon^1\otimes_\tau \alpha_\varepsilon^2
  \in
  \mathfrak A[t,T].
\)
By \eqref{eq:cost-splitting-concatenation-poisson},
\begin{align*}
  \mathcal J_{t,T}(\zeta,\alpha_\varepsilon)
  &=
  \mathcal J_{t,\tau}(\zeta,\alpha_\varepsilon^1)
  +
  \mathcal J_{\tau,T}(\vartheta_\varepsilon,\alpha_\varepsilon^2)
  \\
  &\le
  \mathcal J_{t,\tau}(\zeta,\alpha_\varepsilon^1)
  +
  V(\tau,\vartheta_\varepsilon)
  +
  \varepsilon
  \\
  &\le
  \Xi(t,\tau,\zeta)+2\varepsilon.
\end{align*}
Taking the infimum over \(\alpha\in\mathfrak A[t,T]\), we obtain
\(
  V(t,\zeta)\le \Xi(t,\tau,\zeta)+2\varepsilon.
\)
Since \(\varepsilon>0\) is arbitrary,
\[
  V(t,\zeta)\le \Xi(t,\tau,\zeta).
\]

Combining the two inequalities proves \eqref{eq:DPP-general-poisson}. The
equivalent form \eqref{eq:DPP-general-J-poisson} is immediate from the
definition of \(\mathcal J_{t,\tau}\).
\end{proof}

\section{Linear data and exact moment reduction under mean-dependent affine-Gaussian refresh rules}
\label{sec:LG-moment}

We now specialize the formulation of Section~\ref{sec:dpp-poisson} to linear
state dynamics and to a mean-dependent affine-Gaussian subclass of admissible
Poisson-refresh sampling rules \cite{Wonham68,YongZhou1999StochasticControls}. This subclass remains within the admissible framework of Sections~\ref{sec:poisson-refresh-formulation}--\ref{sec:dpp-poisson}, allowing the general well-posedness, entropy, and dynamic programming results developed there to be specialised to this tractable linear setting.

Because the control is refreshed at random times and then held constant between
refreshes, the full joint law of \((X_t,U_t)\) is not generally Gaussian, even
when the state dynamics are linear and the refresh laws are Gaussian. The
finite-dimensional reduction below is therefore a moment reduction: it tracks
the exact evolution of the first and second moments within the chosen
affine-Gaussian class.

Throughout this section, \(\widehat V\) denotes the value function of this
restricted moment problem. In general, \(\widehat V\) should not be identified
with the unrestricted value function \(V\) of Section~\ref{sec:dpp-poisson}.

We specialize the coefficients to the linear drift
\begin{equation}
  f(x,u)=Ax+Bu,
  \qquad
  A\in\mathbb R^{n\times n},\quad
  B\in\mathbb R^{n\times q},
  \label{eq:LG-primitive-drift-poisson}
\end{equation}
and take the diffusion covariance to be constant,
\begin{equation}
  \Sigma(x,u)\equiv Q,
  \qquad
  Q\in\mathbb S_{++}^n.
  \label{eq:LG-primitive-covariance-poisson}
\end{equation}
Equivalently, the state dynamics are
\begin{equation}
  \mathrm dX_t
  =
  (A X_t+B U_t)\,\mathrm dt
  +
  Q^{1/2}\,\mathrm dW_t.
  \label{eq:LG-state-SDE-poisson}
\end{equation}

The reference drift is taken to be
\begin{equation}
  f^I(t,x,u)
  =
  A^I x+B^I u+x_t^r,
  \qquad
  A^I\in\mathbb R^{n\times n},
  \quad
  B^I\in\mathbb R^{n\times q},
  \label{eq:LG-reference-drift-poisson}
\end{equation}
where \(x^r=(x_t^r)_{t\in[0,T]}\) is continuous. Set
\(\widetilde A:=A-A^I\), \(\widetilde B:=B-B^I\), and \(R:=Q^{-1}\).
Then the local state cost \eqref{eq:local-state-integrand-poisson} becomes
\begin{equation}
  L(t,x,u)
  =
  \frac12
  \bigl(
    \widetilde A x+\widetilde B u-x_t^r
  \bigr)^\top
  R
  \bigl(
    \widetilde A x+\widetilde B u-x_t^r
  \bigr).
  \label{eq:LG-state-cost-poisson}
\end{equation}

The reference refresh kernel is Gaussian and independent of both \(x\) and
\(\mu\):
\begin{equation}
  c_t^I(\cdot\mid x,\mu)\equiv c_t^I(\cdot)
  =
  \mathcal N(u_t^r,\Gamma),
  \qquad
  \Gamma\in\mathbb S_{++}^q,
  \label{eq:LG-reference-refresh-kernel}
\end{equation}
where \(u^r=(u_t^r)_{t\in[0,T]}\) is continuous.

We now define the restricted mean-dependent affine-Gaussian refresh rules used
in the moment reduction. These rules depend on the current state law only
through its mean. Thus, for \(\mu\in\mathcal P_2(\mathbb R^n)\), write
\(
  m_\mu:=\int_{\mathbb R^n}x\,\mu(\mathrm dx)
\).

Fix a Borel map \(g:[0,1]\to\mathbb R^q\) such that the pushforward of
Lebesgue measure on \([0,1]\) by \(g\) is the standard Gaussian law
\(\mathcal N(0,I_q)\). In particular,
\begin{equation}
  \int_0^1 g(z)\,\mathrm dz=0,
  \qquad
  \int_0^1 g(z)g(z)^\top\,\mathrm dz=I_q,
  \qquad
  \int_0^1 \|g(z)\|^2\,\mathrm dz=q.
  \label{eq:fixed-gaussian-realizer}
\end{equation}

We fix structural bounds
\(
  \bar K,\bar a,L_a,\underline c_S,\overline c_S,L_S>0
\)
so that every rule satisfying the conditions below induces a sampling rule
satisfying the admissibility estimates of \(\mathfrak A[t,T]\) with the fixed
constants \(L_\alpha\) and \(M_\alpha\).

A mean-dependent affine-Gaussian refresh rule on \([t,T]\) is specified by three
Borel maps
\[
  K:[t,T]\to\mathbb R^{q\times n},
  \qquad
  a:[t,T]\times\mathbb R^n\to\mathbb R^q,
  \qquad
  S:[t,T]\times\mathbb R^n\to\mathbb S_{++}^q,
\]
satisfying the following conditions. We write
\(
  \beta=(K,a,S)
\)
for the resulting triple. 
\begin{enumerate}
\item
\(
  \sup_{r\in[t,T]}\|K(r)\|\le \bar K
\).

\item
\(a\) is bounded and globally Lipschitz in its second argument, uniformly in
time:
\[
  \|a(r,m)\|\le \bar a,
  \qquad
  \|a(r,m)-a(r,m')\|\le L_a\|m-m'\|,
\]
for all \(r\in[t,T]\) and \(m,m'\in\mathbb R^n\).

\item
\(S\) is bounded and globally Lipschitz in its second argument, uniformly in
time, and satisfies the uniform ellipticity bounds
\begin{equation}
  \underline c_S I_q
  \le
  S(r,m)
  \le
  \overline c_S I_q,
  \qquad
  \|S(r,m)-S(r,m')\|\le L_S\|m-m'\|,
  \label{eq:LG-uniform-S-poisson}
\end{equation}
for all \(r\in[t,T]\) and \(m,m'\in\mathbb R^n\).
\end{enumerate}
The structural bounds are chosen so that the admissibility estimates in
Lemma~\ref{lem:LG-affine-rule-admissible} are controlled by the fixed constants
\(L_\alpha\) and \(M_\alpha\) in
\eqref{eq:sampling-rule-lipschitz}--\eqref{eq:sampling-rule-second-moment}.
The set of such triples is denoted by
\(
  \mathfrak B_{\mathrm{aff}}[t,T]
\).

For \(\beta=(K,a,S)\in\mathfrak B_{\mathrm{aff}}[t,T]\), define
\begin{equation}
  \alpha_s^\beta(x,\mu,z)
  :=
  a(s,m_\mu)
  -
  K(s)(x-m_\mu)
  +
  S(s,m_\mu)^{1/2}g(z),
  \label{eq:LG-affine-Gaussian-sampling-rule}
\end{equation}
for \(s\in[t,T]\), \(x\in\mathbb R^n\),
\(\mu\in\mathcal P_2(\mathbb R^n)\), and \(z\in[0,1]\). The induced
post-refresh kernel is
\begin{equation}
  \gamma_s^\beta(\mathrm du\mid x,\mu)
  =
  \mathcal N\!\bigl(
    a(s,m_\mu)-K(s)(x-m_\mu),
    S(s,m_\mu)
  \bigr).
  \label{eq:LG-affine-Gaussian-refresh-kernel}
\end{equation}

\begin{lemma}\label{lem:LG-affine-rule-admissible}
For every \(\beta\in\mathfrak B_{\mathrm{aff}}[t,T]\), the sampling rule
\(\alpha^\beta\) defined in \eqref{eq:LG-affine-Gaussian-sampling-rule}
belongs to \(\mathfrak A[t,T]\).
\end{lemma}

\begin{proof}
The map \((s,x,\mu,z)\mapsto\alpha_s^\beta(x,\mu,z)\) is Borel measurable. Also,
for \(\mu,\mu'\in\mathcal P_2(\mathbb R^n)\),
\[
  \|m_\mu-m_{\mu'}\|
  \le
  W_2(\mu,\mu').
\]
Since the matrix square-root map is Lipschitz on the set of symmetric matrices
with spectrum contained in \([\underline c_S,\overline c_S]\), there exists a
constant \(C>0\), depending only on the structural bounds defining
\(\mathfrak B_{\mathrm{aff}}[t,T]\), such that
\[
  \|S(s,m)^{1/2}-S(s,m')^{1/2}\|
  \le
  C\|m-m'\|.
\]
Using this estimate together with the boundedness of \(K\), the Lipschitz
property of \(a\), and \eqref{eq:fixed-gaussian-realizer}, we obtain
\[
  \int_0^1
  \left\|
    \alpha_s^\beta(x,\mu,z)-\alpha_s^\beta(x',\mu',z)
  \right\|^2\,\mathrm dz
  \le
  C\Bigl(
    \|x-x'\|^2+W_2(\mu,\mu')^2
  \Bigr).
\]
Since
\[
  \|x-x'\|^2+W_2(\mu,\mu')^2
  \le
  \bigl(\|x-x'\|+W_2(\mu,\mu')\bigr)^2,
\]
and since the compatibility condition ensures that the constant in the
preceding estimate is bounded by \(L_\alpha^2\), we obtain
\[
  \int_0^1
  \left\|
    \alpha_s^\beta(x,\mu,z)-\alpha_s^\beta(x',\mu',z)
  \right\|^2\,\mathrm dz
  \le
  L_\alpha^2
  \Bigl(
    \|x-x'\|+W_2(\mu,\mu')
  \Bigr)^2.
\]

Similarly, using the boundedness of \(K\), \(a\), and \(S\), together with
\eqref{eq:fixed-gaussian-realizer}, we get
\[
  \int_0^1
  \left\|
    \alpha_s^\beta(x,\mu,z)
  \right\|^2\,\mathrm dz
  \le
  C\bigl(
    1+\|x\|^2+\|m_\mu\|^2
  \bigr).
\]
Since \(\|m_\mu\|^2\le m_2(\mu)\), and since the compatibility condition ensures
that the constant in the preceding estimate is bounded by \(M_\alpha\), it
follows that
\[
  \int_0^1
  \left\|
    \alpha_s^\beta(x,\mu,z)
  \right\|^2\,\mathrm dz
  \le
  M_\alpha
  \bigl(
    1+\|x\|^2+m_2(\mu)
  \bigr).
\]
Thus \(\alpha^\beta\) satisfies the measurability, Lipschitz, and second-moment
requirements in the definition of \(\mathfrak A[t,T]\).
\end{proof}

\subsection{Moment state and reduced dynamics}
Having fixed the mean-dependent affine-Gaussian refresh class, we now record the
finite-dimensional variables whose evolution closes under the linear
Poisson-refresh dynamics. For a law
\(
  \zeta\in\mathcal P_2(\mathbb R^n\times\mathbb R^q),
\)
define its joint first and second moments by
\begin{align}
  m_\zeta^x
  &:=
  \int_{\mathbb R^n\times\mathbb R^q}
    x\,\zeta(\mathrm dx,\mathrm du),
  \label{eq:joint-mean-x}
  \\
  m_\zeta^u
  &:=
  \int_{\mathbb R^n\times\mathbb R^q}
    u\,\zeta(\mathrm dx,\mathrm du),
  \label{eq:joint-mean-u}
  \\
  \Sigma_\zeta^{xx}
  &:=
  \int_{\mathbb R^n\times\mathbb R^q}
    (x-m_\zeta^x)(x-m_\zeta^x)^\top\,
    \zeta(\mathrm dx,\mathrm du),
  \label{eq:joint-cov-xx}
  \\
  C_\zeta^{xu}
  &:=
  \int_{\mathbb R^n\times\mathbb R^q}
    (x-m_\zeta^x)(u-m_\zeta^u)^\top\,
    \zeta(\mathrm dx,\mathrm du),
  \label{eq:joint-cov-xu}
  \\
  \Omega_\zeta^{uu}
  &:=
  \int_{\mathbb R^n\times\mathbb R^q}
    (u-m_\zeta^u)(u-m_\zeta^u)^\top\,
    \zeta(\mathrm dx,\mathrm du).
  \label{eq:joint-cov-uu}
\end{align}
We also write
\begin{equation}
  \Xi_\zeta
  :=
  \begin{pmatrix}
    \Sigma_\zeta^{xx} & C_\zeta^{xu} \\
    (C_\zeta^{xu})^\top & \Omega_\zeta^{uu}
  \end{pmatrix}
  \in\mathbb S_+^{n+q}.
  \label{eq:joint-covariance-block}
\end{equation}

The reduced moment state space is
\[
  \mathfrak Z
  :=
  \left\{
    (m,\bar u,\Sigma,C,\Omega):
    m\in\mathbb R^n,\ 
    \bar u\in\mathbb R^q,\ 
    \Sigma\in\mathbb S_+^n,\ 
    C\in\mathbb R^{n\times q},\
    \Omega\in\mathbb S_+^q,\
    \begin{pmatrix}
      \Sigma & C \\
      C^\top & \Omega
    \end{pmatrix}
    \in\mathbb S_+^{n+q}
  \right\}.
\]

For \(\mathfrak z =(m,\bar u,\Sigma,C,\Omega)\in\mathfrak Z\), define
\begin{align}
  F^m(\mathfrak z)
  &:=
  A m+B\bar u,
  \label{eq:LG-reduced-drift-m}
  \\
  F^{\bar u}(\mathfrak z;a)
  &:=
  \lambda(a-\bar u),
  \label{eq:LG-reduced-drift-ubar}
  \\
  F^\Sigma(\mathfrak z)
  &:=
  A\Sigma+\Sigma A^\top+BC^\top+CB^\top+Q,
  \label{eq:LG-reduced-drift-Sigma}
  \\
  F^C(\mathfrak z;k)
  &:=
  AC+B\Omega-\lambda C-\lambda\Sigma k^\top,
  \label{eq:LG-reduced-drift-C}
  \\
  F^\Omega(\mathfrak z;k,a,\Theta)
  &:=
  \lambda\Bigl(
    k\Sigma k^\top
    +
    \Theta
    +
    (a-\bar u)(a-\bar u)^\top
    -
    \Omega
  \Bigr).
  \label{eq:LG-reduced-drift-Omega}
\end{align}

For \(\beta=(K,a,S)\in\mathfrak B_{\mathrm{aff}}[t_0,T]\) and a moment
trajectory
\(
\mathfrak z_s=(m_s,\bar u_s,\Sigma_s,C_s,\Omega_s),
\)
set
\begin{equation}
  K_s:=K(s),
  \qquad
  a_s:=a(s,m_s),
  \qquad
  S_s:=S(s,m_s).
  \label{eq:LG-reduced-pointwise-coefficients}
\end{equation}
The reduced deterministic dynamics are
\begin{align}
  \dot m_s
  &=
  F^m(\mathfrak z_s),
  \label{eq:LG-moment-ODE-m}
  \\
  \dot{\bar u}_s
  &=
  F^{\bar u}(\mathfrak z_s;a_s),
  \label{eq:LG-moment-ODE-ubar}
  \\
  \dot\Sigma_s
  &=
  F^\Sigma(\mathfrak z_s),
  \label{eq:LG-moment-ODE-Sigma}
  \\
  \dot C_s
  &=
  F^C(\mathfrak z_s;K_s),
  \label{eq:LG-moment-ODE-C}
  \\
  \dot\Omega_s
  &=
  F^\Omega(\mathfrak z_s;K_s,a_s,S_s).
  \label{eq:LG-moment-ODE-Omega}
\end{align}

The next proposition shows that the deterministic moment system is not a
closure approximation: it is exactly realised by the corresponding
Poisson-refresh jump--diffusion for every initial law with the prescribed first
and second moments.

\begin{proposition}[Reduced ODE well-posedness and exact realization]
\label{prop:LG-moment-closure}
Fix \(t_0\in[0,T]\), let
\(
  \mathfrak z^0=(m^0,\bar u^0,\Sigma^0,C^0,\Omega^0)\in\mathfrak Z,
\)
and let
\(
  \beta=(K,a,S)\in\mathfrak B_{\mathrm{aff}}[t_0,T].
\)
Then the system
\eqref{eq:LG-moment-ODE-m}--\eqref{eq:LG-moment-ODE-Omega}, started from
\(\mathfrak z^0\), admits a unique absolutely continuous solution
\[
  \mathfrak z^{t_0,\mathfrak z^0,\beta}
  =
  (m_s,\bar u_s,\Sigma_s,C_s,\Omega_s)_{s\in[t_0,T]}
  \in C([t_0,T];\mathfrak Z).
\]

Moreover, if \((X_{t_0},U_{t_0})\) is any square-integrable random vector with
first and second moments \(\mathfrak z^0\), independent of an \(n\)-dimensional Brownian
motion \(W\) and of a Poisson random measure \(N(\mathrm ds,\mathrm d\xi)\) on
\([t_0,T]\times[0,1]\) with intensity
\(
  \lambda\,\mathrm ds\,\mathrm d\xi,
\)
then the linear jump--diffusion system
\begin{align}
  \mathrm dX_s
  &=
  (A X_s+B U_s)\,\mathrm ds
  +
  Q^{1/2}\,\mathrm dW_s,
  \label{eq:LG-realized-state-equation}
  \\
  U_s
  &=
  U_{t_0}
  +
  \int_{(t_0,s]\times[0,1]}
    \Bigl(
      a_r-K_r(X_{r-}-m_r)+S_r^{1/2}g(\xi)-U_{r-}
    \Bigr)\,
    N(\mathrm dr,\mathrm d\xi),
  \qquad s\in[t_0,T],
  \label{eq:LG-realized-control-equation}
\end{align}
admits a unique strong solution with finite second moments, and its joint first
and second moments are exactly
\[
  \mathbb E[X_s]=m_s,
  \qquad
  \mathbb E[U_s]=\bar u_s,
\]
\[
  \operatorname{Cov}(X_s)=\Sigma_s,
  \qquad
  \operatorname{Cov}(X_s,U_s)=C_s,
  \qquad
  \operatorname{Cov}(U_s)=\Omega_s,
  \qquad
  s\in[t_0,T].
\]
In particular, the reduced trajectory depends only on the initial first and
second moments and on \(\beta\), not on higher moments of the initial law.
\end{proposition}

\begin{proof}
The coefficients in the reduced ODE are Borel in time and locally Lipschitz in
\((m,\bar u,\Sigma,C,\Omega)\). Hence the system has a unique maximal
absolutely continuous solution on some interval \([t_0,\tau_{\max})\).

Along this local solution, the coefficients \(K_s\), \(a_s\), and \(S_s\) are
bounded Borel functions of time, with \(S_s\) uniformly positive definite. The
linear jump--diffusion
\eqref{eq:LG-realized-state-equation}--%
\eqref{eq:LG-realized-control-equation}
therefore has a unique strong solution on compact subintervals of
\([t_0,\tau_{\max})\).

Applying It\^o's formula for jump--diffusions to the coordinate functions
\(x\) and \(u\), and using \eqref{eq:fixed-gaussian-realizer}, gives
\[
  \frac{\mathrm d}{\mathrm ds}\mathbb E[X_s]
  =
  A\mathbb E[X_s]+B\mathbb E[U_s],
\]
and
\[
  \frac{\mathrm d}{\mathrm ds}\mathbb E[U_s]
  =
  \lambda\Bigl(
    a_s-K_s(\mathbb E[X_s]-m_s)-\mathbb E[U_s]
  \Bigr).
\]
Comparing with the first two ODEs in
\eqref{eq:LG-moment-ODE-m}--\eqref{eq:LG-moment-ODE-ubar} and using uniqueness
for the resulting linear difference system gives \(\mathbb E[X_s]=m_s,\) and \(\mathbb E[U_s]=\bar u_s.\)

Applying It\^o's formula to \(xx^\top\), \(xu^\top\), and \(uu^\top\), and then
subtracting the derivatives of the corresponding products of means, gives
\begin{align*}
  \frac{\mathrm d}{\mathrm ds}\operatorname{Cov}(X_s)
  &=
  A\operatorname{Cov}(X_s)
  +
  \operatorname{Cov}(X_s)A^\top
  +
  B\operatorname{Cov}(X_s,U_s)^\top
  +
  \operatorname{Cov}(X_s,U_s)B^\top
  +
  Q,
  \\
  \frac{\mathrm d}{\mathrm ds}\operatorname{Cov}(X_s,U_s)
  &=
  A\operatorname{Cov}(X_s,U_s)
  +
  B\operatorname{Cov}(U_s)
  -
  \lambda\operatorname{Cov}(X_s,U_s)
  -
  \lambda\operatorname{Cov}(X_s)K_s^\top,
  \\
  \frac{\mathrm d}{\mathrm ds}\operatorname{Cov}(U_s)
  &=
  \lambda\Bigl(
    K_s\operatorname{Cov}(X_s)K_s^\top
    +
    S_s
    +
    (a_s-\bar u_s)(a_s-\bar u_s)^\top
    -
    \operatorname{Cov}(U_s)
  \Bigr).
\end{align*}
These are exactly the last three equations in
\eqref{eq:LG-moment-ODE-Sigma}--\eqref{eq:LG-moment-ODE-Omega}. Hence the
moment identities hold on the local interval.

It remains to prove global existence. Define
\[
  M_s
  :=
  1+\|m_s\|^2+\|\bar u_s\|^2+\operatorname{tr}(\Sigma_s)+\operatorname{tr}(\Omega_s).
\]
Since the block covariance matrix is positive semidefinite,
\[
  \|C_s\|_{\mathrm F}^2
  \le
  \operatorname{tr}(\Sigma_s)\operatorname{tr}(\Omega_s).
\]
Using the ODE system and the boundedness of \(K\), \(a\), and \(S\), one obtains
\[
  \frac{\mathrm d}{\mathrm ds}M_s
  \le
  C M_s
\]
for a constant \(C\) depending only on the linear data, the refresh rate, and
the structural bounds defining \(\mathfrak B_{\mathrm{aff}}[t_0,T]\). Thus
\(M_s\) remains bounded on finite intervals, and the maximal solution cannot
explode before \(T\). Therefore \(\tau_{\max}=T\). The exact realization of the
moment trajectory follows from the moment identities above, while the strong
uniqueness of the jump--diffusion follows from the standard well-posedness
theory for linear jump--diffusions with globally Lipschitz coefficients.
\end{proof}

\subsection{Running cost on moment space}

The moment dynamics in Proposition~\ref{prop:LG-moment-closure} close on
\(\mathfrak Z\). We now record that, within the same affine-Gaussian refresh
class, the local running cost also depends only on the moment state and the
current refresh parameters.

For a moment state
\(
  \mathfrak z=(m,\bar u,\Sigma,C,\Omega)\in\mathfrak Z,
\)
define
\[
  \Xi(z):=
  \begin{pmatrix}
    \Sigma & C \\
    C^\top & \Omega
  \end{pmatrix},
  \qquad
  \bar y(z):=
  \binom{m}{\bar u},
  \qquad
  G:=
  \begin{pmatrix}
    \widetilde A & \widetilde B
  \end{pmatrix}.
\]

\begin{proposition}[Running cost on moment space]
\label{prop:LG-running-cost-moments}
Fix \(t\in[0,T]\). Let
\(\mathfrak z=(m,\bar u,\Sigma,C,\Omega)\in\mathfrak Z\), and let
\(
  k\in\mathbb R^{q\times n}\), 
\(  a\in\mathbb R^q\), and 
 \( S\in\mathbb S_{++}^q \).
Then the local control-side cost associated with the refresh kernel
\(
  x\longmapsto \mathcal N(a-k(x-m),S)
\)
is
\begin{align}
  \ell_t^{\mathrm{ctrl}}(m,\Sigma;k,a,S)
  &=
  \frac{\lambda}{2}\Bigl(
    \operatorname{tr}(\Gamma^{-1}S)
    -\log\det S
    +\log\det\Gamma
    -q
  \Bigr)
  \notag\\
  &\quad+
  \frac{\lambda}{2}
  \operatorname{tr}\!\bigl(
    k^\top\Gamma^{-1}k\,\Sigma
  \bigr)
  +
  \frac{\lambda}{2}
  (a-u_t^r)^\top\Gamma^{-1}(a-u_t^r).
  \label{eq:LG-running-cost-control-moments}
\end{align}
The local state-side cost is
\begin{equation}
  \ell_t^{\mathrm{state}}(\mathfrak z)
  =
  \frac12
  \operatorname{tr}\!\bigl(
    G^\top R G\,\Xi(\mathfrak z)
  \bigr)
  +
  \frac12
  \bigl(
    G\bar y(\mathfrak z)-x_t^r
  \bigr)^\top
  R
  \bigl(
    G\bar y(\mathfrak z)-x_t^r
  \bigr).
  \label{eq:LG-running-cost-state-moments}
\end{equation}
The total running cost is
\begin{equation}
  \ell_t(\mathfrak z;k,a,S)
  :=
  \ell_t^{\mathrm{ctrl}}(m,\Sigma;k,a,S)
  +
  \ell_t^{\mathrm{state}}(\mathfrak z).
  \label{eq:LG-running-cost-total-moments}
\end{equation}
\end{proposition}

\begin{proof}
Let \((X,U)\) be any square-integrable random vector with moment state
\(\mathfrak z=(m,\bar u,\Sigma,C,\Omega)\). All expectations in this proof are
taken with respect to the law of \((X,U)\).

The control-side formula follows from the Gaussian relative-entropy formula
applied, for fixed \(x\), to
\(
  \mathcal N(a-k(x-m),S) \) and
\(  \mathcal N(u_t^r,\Gamma)
\),
and then averaged over \(X\). The cross term vanishes because \(\mathbb E[X-m]=0\), and
\[
  \mathbb E\!\left[
    (X-m)^\top k^\top\Gamma^{-1}k(X-m)
  \right]
  =
  \operatorname{tr}\!\bigl(k^\top\Gamma^{-1}k\,\Sigma\bigr).
\]

For the state-side term, let \(Y=(X,U)^\top\) have mean \(\bar y(\mathfrak z)\) and
covariance \(\Xi(\mathfrak z)\). Since
\(
  \widetilde A X+\widetilde B U-x_t^r
  =
  GY-x_t^r,
\)
the standard second-moment identity
\[
  \mathbb E\!\left[
    (GY-x_t^r)^\top R(GY-x_t^r)
  \right]
  =
  \operatorname{tr}\!\bigl(G^\top R G\,\Xi(\mathfrak z)\bigr)
  +
  \bigl(G\bar y(\mathfrak z)-x_t^r\bigr)^\top
  R
  \bigl(G\bar y(\mathfrak z)-x_t^r\bigr)
\]
gives \eqref{eq:LG-running-cost-state-moments}.
\end{proof}

\subsection{Reduced deterministic control problem on moment space}

Using the moment dynamics of Proposition~\ref{prop:LG-moment-closure} and the
moment-space running cost of Proposition~\ref{prop:LG-running-cost-moments}, we
now define the restricted deterministic control problem on \(\mathfrak Z\).

For \((t,\mathfrak z)\in[0,T]\times\mathfrak Z\) and
\(\beta=(K,a,S)\in\mathfrak B_{\mathrm{aff}}[t,T]\), let
\(\mathfrak z^{t,\mathfrak z,\beta}=(\mathfrak z_s^{t,\mathfrak z,\beta})_{s\in[t,T]}\), with
\(\mathfrak z_s^{t,\mathfrak z,\beta}=(m_s,\bar u_s,\Sigma_s,C_s,\Omega_s)\), denote the unique
solution of \eqref{eq:LG-moment-ODE-m}--\eqref{eq:LG-moment-ODE-Omega} started
from \(\mathfrak z\). Along this trajectory, write \(K_s:=K(s)\),
\(a_s:=a(s,m_s)\), and \(S_s:=S(s,m_s)\).
The reduced finite-horizon cost is defined by integrating the running cost
\(\ell\) from Proposition~\ref{prop:LG-running-cost-moments} along this moment
trajectory:
\begin{equation}
  \widehat{\mathcal J}_{t,T}(\mathfrak z,\beta)
  :=
  \int_t^T
    \ell_s\!\bigl(
      \mathfrak z_s^{t,\mathfrak z,\beta};K_s,a_s,S_s
    \bigr)\,
    \mathrm ds.
  \label{eq:LG-reduced-cost-moments}
\end{equation}
The reduced value function is
\begin{equation}
  \widehat V(t,\mathfrak z)
  :=
  \inf_{\beta\in\mathfrak B_{\mathrm{aff}}[t,T]}
  \widehat{\mathcal J}_{t,T}(\mathfrak z,\beta),
  \qquad
  (t,\mathfrak z)\in[0,T]\times\mathfrak Z,
  \label{eq:LG-reduced-value-moments}
\end{equation}
with terminal condition \(\widehat V(T,\mathfrak z)=0\) for all \(\mathfrak z \in\mathfrak Z\).

For \(0\le t\le \tau\le T\) and
\(\beta\in\mathfrak B_{\mathrm{aff}}[t,\tau]\), we similarly define the
cost accumulated over the subinterval \([t,\tau]\) by
\begin{equation}
  \widehat{\mathcal J}_{t,\tau}(\mathfrak z,\beta)
  :=
  \int_t^\tau
    \ell_s\!\bigl(
      \mathfrak z_s^{t,\mathfrak z,\beta};K_s,a_s,S_s
    \bigr)\,\mathrm ds.
  \label{eq:LG-reduced-short-cost-moments}
\end{equation}

\begin{proposition}[Reduced dynamic programming principle]
\label{prop:LG-reduced-DPP-moments}
For every \(0\le t\le\tau\le T\) and every \(\mathfrak z \in\mathfrak Z\),
\begin{equation}
  \widehat V(t,\mathfrak z)
  =
  \inf_{\beta\in\mathfrak B_{\mathrm{aff}}[t,\tau]}
  \left\{
    \widehat{\mathcal J}_{t,\tau}(\mathfrak z,\beta)
    +
    \widehat V\!\bigl(\tau,\mathfrak z_\tau^{t,\mathfrak z,\beta}\bigr)
  \right\}.
  \label{eq:LG-reduced-DPP-moments}
\end{equation}
\end{proposition}

\begin{proof}
The proof is the finite-dimensional deterministic analogue of
Theorem~\ref{thm:DPP-poisson}. Restricting and concatenating triples
\(\beta=(K,a,S)\) preserves membership in \(\mathfrak B_{\mathrm{aff}}\), and
the reduced dynamics are deterministic. Splitting the integral at \(\tau\) and
using the usual \(\varepsilon\)-optimal concatenation argument gives
\eqref{eq:LG-reduced-DPP-moments}.
\end{proof}

\subsection{Reduced Hamiltonian and verification on moment space}
The reduced dynamic programming principle leads to the following verification
formulation of the associated Hamiltonian problem on moment space, in the spirit
of classical stochastic-control verification and HJB arguments
\cite{FlemingSoner2006Controlled}. The goal is to give conditions under which a
smooth candidate function bounds, and in the attained case equals, the reduced
value function \(\widehat V\).

Let
\(
  \widehat W:[0,T]\times
  \mathbb R^n\times\mathbb R^q\times\mathbb S^n\times
  \mathbb R^{n\times q}\times\mathbb S^q
  \to\mathbb R
\)
be smooth. Write \(\nabla_m\widehat W\) and
\(\nabla_{\bar u}\widehat W\) for the gradients with respect to \(m\) and
\(\bar u\), and define \(\nabla_\Sigma\widehat W\),
\(\nabla_C\widehat W\), and \(\nabla_\Omega\widehat W\) through the trace
pairings in the matrix variables.

For pointwise values
\(k\in\mathbb R^{q\times n}\), \(a\in\mathbb R^q\), and
\(\Theta\in\mathbb S_{++}^q\), define the first-order operator associated with
the reduced moment dynamics by
\begin{align}
  \widehat{\mathcal L}^{k,a,\Theta}\widehat W(t,\mathfrak z)
  &:=
  \nabla_m\widehat W(t,\mathfrak z)^\top F^m(\mathfrak z)
  +
  \nabla_{\bar u}\widehat W(t,\mathfrak z)^\top F^{\bar u}(\mathfrak z;a)
  \notag\\
  &\quad+
  \operatorname{tr}\!\bigl(
    \nabla_\Sigma\widehat W(t,\mathfrak z)\,F^\Sigma(\mathfrak z)
  \bigr)
  +
  \operatorname{tr}\!\bigl(
    \nabla_C\widehat W(t,\mathfrak z)^\top F^C(\mathfrak z;k)
  \bigr)
  \notag\\
  &\quad+
  \operatorname{tr}\!\bigl(
    \nabla_\Omega\widehat W(t,\mathfrak z)\,F^\Omega(\mathfrak z;k,a,\Theta)
  \bigr).
  \label{eq:LG-reduced-generator-moments}
\end{align}

Set
\(
  \mathcal A_{\mathrm{pt}}
  :=
  \left\{
    (k,a,\Theta):
    \|k\|\le \bar K,\ 
    \|a\|\le \bar a,\ 
    \underline c_S I_q\le \Theta\le \overline c_S I_q
  \right\}
\).
The corresponding reduced pointwise Hamiltonian is
\begin{align}
  \widehat H(t,\mathfrak z,p_m,p_{\bar u},P_\Sigma,P_C,P_\Omega)
  &:=
  \inf_{(k,a,\Theta)\in\mathcal A_{\mathrm{pt}}}
  \Bigl\{
    p_m^\top F^m(\mathfrak z)
    +
    p_{\bar u}^\top F^{\bar u}(\mathfrak z;a)
    +
    \operatorname{tr}\!\bigl(P_\Sigma F^\Sigma(\mathfrak z)\bigr)
  \notag\\
  &\qquad+
    \operatorname{tr}\!\bigl(P_C^\top F^C(\mathfrak z;k)\bigr)
    +
    \operatorname{tr}\!\bigl(P_\Omega F^\Omega(\mathfrak z;k,a,\Theta)\bigr)
    +
    \ell_t(\mathfrak z;k,a,\Theta)
  \Bigr\}.
  \label{eq:LG-reduced-Hamiltonian-moments}
\end{align}

Whenever the Hamiltonian is finite, the associated reduced Hamilton--Jacobi
equation is
\begin{equation}
  0
  =
  \partial_t\widehat W(t,\mathfrak z)
  +
  \widehat H\!\bigl(
    t,\mathfrak z,
    \nabla_m\widehat W(t,\mathfrak z),
    \nabla_{\bar u}\widehat W(t,\mathfrak z),
    \nabla_\Sigma\widehat W(t,\mathfrak z),
    \nabla_C\widehat W(t,\mathfrak z),
    \nabla_\Omega\widehat W(t,\mathfrak z)
  \bigr),
  \label{eq:LG-reduced-HJB-moments}
\end{equation}
with terminal condition \(\widehat W(T,\mathfrak z)=0\).

\begin{theorem}[Reduced verification on moment space]
\label{thm:LG-reduced-verification-moments}
Let \(\widehat W\) be continuous, continuously differentiable in time, and
continuously differentiable in the moment variables. Assume that
\(\widehat W(T,\mathfrak z)=0\) for all \(\mathfrak z\in\mathfrak Z\) and that
\begin{equation}
  \partial_t\widehat W(t,\mathfrak z)
  +
  \widehat H\!\bigl(
    t,\mathfrak z,
    \nabla_m\widehat W(t,\mathfrak z),
    \nabla_{\bar u}\widehat W(t,\mathfrak z),
    \nabla_\Sigma\widehat W(t,\mathfrak z),
    \nabla_C\widehat W(t,\mathfrak z),
    \nabla_\Omega\widehat W(t,\mathfrak z)
  \bigr)
  \ge 0
  \label{eq:LG-reduced-HJB-super-moments}
\end{equation}
for all \((t,\mathfrak z)\in[0,T)\times\mathfrak Z\). Then
\begin{equation}
  \widehat W(t,\mathfrak z)\le \widehat V(t,\mathfrak z),
  \qquad
  (t,\mathfrak z)\in[0,T]\times\mathfrak Z.
  \label{eq:LG-reduced-verification-ineq-moments}
\end{equation}
If, in addition, equality holds in
\eqref{eq:LG-reduced-HJB-super-moments} and there exists
\(
  \widehat\beta=(\widehat K,\widehat a,\widehat S)
  \in\mathfrak B_{\mathrm{aff}}[0,T]
\)
such that, for each initial pair \((t,\mathfrak z)\), the pointwise values
\(
  \widehat K(s) \),
\(  \widehat a(s,\widehat m_s) \), and \(\widehat S(s,\widehat m_s)
\)
attain the infimum in \eqref{eq:LG-reduced-Hamiltonian-moments} for almost
every \(s\in[t,T]\) along the reduced trajectory
\(
  \widehat{\mathfrak z}_s^{t,\mathfrak z}
  =
  (\widehat m_s,\widehat{\bar u}_s,\widehat\Sigma_s,\widehat C_s,\widehat\Omega_s)
\)
generated by \(\widehat\beta\), then
\begin{equation}
  \widehat W(t,\mathfrak z)=\widehat V(t,\mathfrak z),
  \qquad
  (t,\mathfrak z)\in[0,T]\times\mathfrak Z,
  \label{eq:LG-reduced-verification-eq-moments}
\end{equation}
and, for each initial pair \((t,\mathfrak z)\), the rule
\(\widehat\beta|_{[t,T]}\) is optimal for the reduced moment problem started
from \((t,\mathfrak z)\).
\end{theorem}

\begin{proof}
The proof is the standard verification argument for deterministic control. Fix
\((t,\mathfrak z)\in[0,T]\times\mathfrak Z\) and
\(\beta=(K,a,S)\in\mathfrak B_{\mathrm{aff}}[t,T]\). Along the corresponding
trajectory \(\mathfrak z_s^{t,\mathfrak z,\beta}\), write
\(K_s:=K(s) \),  \(a_s:=a(s,m_s)\), and, \(S_s:=S(s,m_s)\).
By the chain rule,
\[
  \frac{\mathrm d}{\mathrm ds}
  \widehat W\!\bigl(s,\mathfrak z_s^{t,\mathfrak z,\beta}\bigr)
  =
  \partial_s\widehat W\!\bigl(s,\mathfrak z_s^{t,\mathfrak z,\beta}\bigr)
  +
  \widehat{\mathcal L}^{K_s,a_s,S_s}
  \widehat W\!\bigl(s,\mathfrak z_s^{t,\mathfrak z,\beta}\bigr)
\]
for almost every \(s\in[t,T]\). Since the Hamiltonian is the infimum over the
admissible pointwise values, \eqref{eq:LG-reduced-HJB-super-moments} gives
\[
  \partial_s\widehat W\!\bigl(s,\mathfrak z_s^{t,\mathfrak z,\beta}\bigr)
  +
  \widehat{\mathcal L}^{K_s,a_s,S_s}
  \widehat W\!\bigl(s,\mathfrak z_s^{t,\mathfrak z,\beta}\bigr)
  +
  \ell_s\!\bigl(
    \mathfrak z_s^{t,\mathfrak z,\beta};K_s,a_s,S_s
  \bigr)
  \ge 0.
\]
Integrating from \(t\) to \(T\) and using
\(\widehat W(T,\cdot)=0\) yields
\[
  \widehat W(t,\mathfrak z)
  \le
  \int_t^T
    \ell_s\!\bigl(
      \mathfrak z_s^{t,\mathfrak z,\beta};K_s,a_s,S_s
    \bigr)\,\mathrm ds.
\]
Taking the infimum over
\(\beta\in\mathfrak B_{\mathrm{aff}}[t,T]\) proves
\eqref{eq:LG-reduced-verification-ineq-moments}.

If the Hamiltonian equality is attained by \(\widehat\beta\), the same
calculation along the trajectory generated by \(\widehat\beta\) gives equality:
\[
  \widehat W(t,\mathfrak z)
  =
  \int_t^T
    \ell_s\!\bigl(
      \widehat{\mathfrak z}_s^{t,\mathfrak z};
      \widehat K(s),
      \widehat a(s,\widehat m_s),
      \widehat S(s,\widehat m_s)
    \bigr)\,\mathrm ds.
\]
Therefore \(\widehat V(t,\mathfrak z)\le \widehat W(t,\mathfrak z)\). Combined
with \eqref{eq:LG-reduced-verification-ineq-moments}, this proves
\eqref{eq:LG-reduced-verification-eq-moments}, and optimality of
\(\widehat\beta\) follows.
\end{proof}

\section{Conclusion}
\label{sec:conc_Rem}

This paper developed a continuous-time control formulation for stochastic
systems with law-dependent Poisson-refresh randomised feedback. The control was
modelled as a sampled process \(U\), held constant between refresh times and
resampled at Poisson times from a feedback rule depending on the current state
and state law. This led to a closed-loop jump--diffusion for the joint process
\((X,U)\), where the state dynamics are driven by the realised control value
rather than by an averaged control kernel.

Under the stated Lipschitz, moment, and ellipticity assumptions, we proved
well-posedness of the closed-loop system for the admissible class of sampling
rules. We also showed that the joint law flow is c\`adl\`ag, the state-law flow
is continuous, the solution law is stable with respect to the initial joint law,
and the admissible class is stable under restriction and concatenation. These
properties provided the structure needed to formulate a continuation problem
with the joint law of \((X_t,U_t)\) as the state variable.

The finite-horizon objective was built from two local contributions. The
control-side term measures the relative entropy of the post-refresh kernel
against a prescribed reference refresh rule. The state-side term measures drift
mismatch relative to a prescribed reference drift and is evaluated at the
realised control value. Combining these terms gave an additive cost functional
on the closed-loop law, and we proved a dynamic programming principle for the
resulting value function on \( \mathcal P_2(\mathbb R^n\times\mathbb R^q) \).

We then demonstrated the formulation on a linear-data example. Under the
admissible mean-dependent affine-Gaussian refresh rules considered there, the
full joint law of \((X_t,U_t)\) need not remain Gaussian, but its first and
second moments close exactly. This yielded a deterministic reduced control
problem on moment space, with an explicit running cost, a reduced dynamic
programming principle, a reduced Hamiltonian, and a verification theorem.

Several questions remain open. One is to compare the Poisson-refresh
formulation more systematically with relaxed-control and path-space entropy
formulations for control-dependent diffusions. Another is to identify further
classes of refresh rules or law families for which finite-dimensional
reductions are exact or analytically tractable.

The dynamic programming principle also points toward a Hamilton--Jacobi theory
on \( \mathcal P_2(\mathbb R^n\times\mathbb R^q) \),
with the joint law of \((X_t,U_t)\) as the continuation variable. Developing a
rigorous full-space Hamilton--Jacobi characterization, including an appropriate
test-function framework and uniqueness theory, is beyond the scope of the
present paper and is left for future work.

\bibliographystyle{amsplain}  
\bibliography{Paper_Ref}

\section*{Declarations}

\begin{itemize}
\item Funding: This work was supported by the Engineering and Physical Sciences Research Council (EPSRC), UK Research and Innovation (UKRI), award reference UKRI2404.
\item Acknowledgements: The author thanks Dr. Roger Tribe for helpful feedback on some aspects on an earlier draft of this paper.
\item Competing Interests: The author have no relevant financial or non-financial interests to disclose.
\item Author Contributions: Randa Herzallah conceived the study, developed the theory and proofs, and wrote the manuscript.
\end{itemize}

\end{document}